\documentclass[onefignum,onetabnum]{siamart190516}

\usepackage{xcolor}
\usepackage{amsfonts}
\usepackage{graphicx}
\usepackage{epstopdf}
\usepackage{overpic}
\usepackage{amsmath,esint}
\usepackage{upquote,enumitem}
\usepackage{mathtools}
\usepackage{listings}
\usepackage{algorithm}
\usepackage{algpseudocode}

\crefname{lstlisting}{listing}{listings}
\Crefname{lstlisting}{Listing}{Listings}

\ifpdf
  \DeclareGraphicsExtensions{.eps,.pdf,.png,.jpg}
\else
  \DeclareGraphicsExtensions{.eps}
\fi

\newsiamremark{remark}{Remark}
\newsiamremark{hypothesis}{Hypothesis}
\crefname{hypothesis}{Hypothesis}{Hypotheses}
\newsiamthm{claim}{Claim}

\headers{Control of bifurcation structures using shape optimization}{Nicolas Boull\'e, Patrick E. Farrell, and Alberto Paganini}

\title{Control of bifurcation structures using shape optimization \thanks{Submitted to the editors \today.
\funding{NB was supported by the EPSRC Centre for Doctoral Training in Industrially Focused Mathematical Modelling (EP/L015803/1) in collaboration with Simula Research Laboratory. PEF was supported by EPSRC grants EP/V001493/1 and EP/R029423/1.}}}

\author{Nicolas Boull\'e\thanks{Mathematical Institute, University of Oxford, Oxford, OX2 6GG, UK. (\email{boulle@maths.ox.ac.uk})}
\and Patrick E. Farrell\thanks{Mathematical Institute, University of Oxford, Oxford, OX2 6GG, UK. (\email{patrick.farrell@maths.ox.ac.uk})}
\and Alberto Paganini\thanks{School of Mathematics and Actuarial Science, University of Leicester, Leicester, LE1 7RH, UK. (\email{a.paganini@leicester.ac.uk})}}

\usepackage{amsopn}

\DeclareMathOperator{\tr}{tr}

\begin{document}

\maketitle

\begin{abstract}
Many problems in engineering can be understood as controlling the bifurcation structure of a given device. For example, one may wish to delay the onset of instability, or bring forward a bifurcation to enable rapid switching between states. We propose a numerical technique for controlling the bifurcation diagram of a nonlinear partial differential equation by varying the shape of the domain. Specifically, we are able to delay or advance a given branch point to a target parameter value. The algorithm consists of solving a shape optimization problem constrained by an augmented system of equations, the Moore--Spence system, that characterize the location of the branch points. Numerical experiments on the Allen--Cahn, Navier--Stokes, and hyperelasticity equations demonstrate the effectiveness of this technique in a wide range of settings.
\end{abstract}



\section{Introduction}\label{sec_intro}

A striking property of nonlinear partial differential equations (PDEs) is that they may support
multiple solutions, not trivially related via the nullspace of a linear
operator. These solutions are captured in a bifurcation diagram,
recording the solutions $u \in U$ to an equation 
\begin{equation}
\label{eq:intro}
F(u, \lambda) = 0
\end{equation}
as
a parameter $\lambda \in \mathbb{R}$ is varied. Here, $U$ is a Banach
space of functions defined on a bounded domain $\Omega \subset \mathbb{R}^d$,  $d \in \mathbb{N}_+$,
and $F: U \times \mathbb{R} \to U^{\star}$, where $U^{*}$ denotes the dual space of $U$\footnote{This is a convenient way of writing variational
problems. The variational problem
\begin{equation*}
\text{find } u \in U \text{ such that } G(u, \lambda; v) = 0 \text{ for all } v \in U
\end{equation*}
may be recast as \cref{eq:intro} by defining
\begin{equation*}
\langle F(u, \lambda), v \rangle \coloneqq G(u, \lambda; v).
\end{equation*}
}. Throughout the paper, we assume that the nonlinear operator $F$ is Fr\'echet differentiable, i.e.~$F \in C^1(U \times \mathbb{R}, U^*)$. When a system is to be
designed, consideration of multiple solutions can be of primary
importance: for example, one may wish to ensure that an aircraft remains
in a high-lift rather than low-lift regime in
takeoff~\cite{kamenetskiy2014}, or may wish to exploit bistability in
the design of a liquid crystal display that only consumes power when an
image is changed~\cite{jones2012}.

We consider the mathematical \emph{control} of the bifurcation structure of a physical system and formulate an algorithm that modifies the domain $\Omega$ on which a PDE is posed in order to advance or delay a single specified bifurcation point. Several engineering problems may be understood as such a task. For instance, in wing design one typically aims to avoid instabilities due to the buckling of the wing structure. In this example, the bifurcation parameter $\lambda$ is the scaling of the load applied to the wing. To avoid buckling, one can modify the design of the structure $\Omega$ to delay the first branch point (with respect to $\lambda$), and thus remove potentially unstable solutions. A similar approach might be taken to the design of aircraft stiffeners~\cite{xia2020nonlinear} and other components. On the contrary, in other applications, inducing new branches might be desirable. For instance, this is the case in the design of new mechanical metamaterials that exploit instabilities to rapidly switch between configurations, as in the snapping of a Venus flytrap~\cite{forterre2005,medina2020navigating,vidoli2008}.

Both bifurcation analysis and PDE-constrained optimization are now mature fields, but their combination has essentially been explored, in a number of related works, through indirect methods that exploit prior model information in specific cases. As an example, Penzler et al.~\cite{penzler2012phase} performed shape optimization in the context of nonlinear elasticity and had to design a robust minimization method to handle non-uniqueness of the solution arising from buckling instabilities. More recently, Thomsen et al.~\cite{thomsen2018buckling} considered a topology optimization problem to maximize strength of materials against buckling instabilities, while the control of flow bifurcations in channels using shape optimization was analyzed in~\cite{wang2019enhanced}. In both cases, the state equations and corresponding eigenvalue problems for finding the critical modes are decoupled in the inner loop of the optimization algorithm to reduce the implementation complexity. Finally, a direct method for performing shape optimization of buckling structures using a characterization of critical points was proposed in~\cite{reitinger1994shape}. The algorithm optimizes the critical load factor over some shape parameters variables using SQP methods~\cite[Chapt.~18]{nocedal2006numerical}, and the underlying system of equations is solved using a penalty formulation to avoid convergence issues near critical points.

In this work, we formulate the task of controlling an isolated branch point as a PDE-constrained shape optimization problem. The essential idea of the formulation is to characterize the branch point as the solution of an augmented system of equations. There are many choices for how to describe branch points in this way (see e.g.~\cite[\S 5.4]{seydel2009practical}). In our case we choose the system proposed by Seydel~\cite{seydel1979} and Moore \& Spence~\cite{moore1980calculation}, referred to as the Moore--Spence system. We then propose a generic numerical algorithm for solving the constrained optimization problem. This is a rather challenging task: none of the state equation \cref{eq:intro}, the Moore--Spence system, or the outer shape optimization problem have unique solutions.

Much of the literature on PDE-constrained optimization requires that the control-to-state
map is single-valued~\cite{Al07,DeZo11,HaMae03,hinze2009,SoZo92,troltzsch2010}. In fact, in many analyses some smoothness
of the map is further required. However, the case where the implicit
function theorem cannot be applied has also been considered, in both finite and infinite dimensions~\cite{dennis1997,heinkenschloss2014,heinkenschloss2002}.
Some experiments involving multiple solutions of optimal control problems, arising
from multiple solutions of the underlying control-to-state map, were recently
presented in \cite{pichi2020driving}.
In our work, we address the challenges posed by the co-existence of multiple
solutions by careful selection of the initial guesses employed and their updates, which
occur while the optimization is carried out, to exploit local uniqueness of isolated branch points.
Our algorithm proves successful in several numerical experiments, and allows us to find a new shape
$\Omega^\star$ for which a branch of solutions arises at a specific
value $\lambda = \lambda^\star$. Additionally, we release an open-source implementation of the algorithm, built on Firedrake~\cite{rathgeber2016firedrake}, PETSc~\cite{balay2020petsc}, and Fireshape~\cite{paganini2020fireshape}.

While this paper focuses on controlling simple pitchfork and fold bifurcation
points with respect to the shape of the domain, we expect that the
ideas developed here also apply to other settings, such as the control
of the bifurcation diagram with respect to a parameter $\lambda_1 \in \mathbb{R}$ as
another parameter $\lambda_2$ in some (possibly infinite-dimensional) parameter
space is varied, or the control of other kinds
of bifurcations such as Hopf bifurcations.

The paper is organized as follows. We first review the characterization we use of
branch points using augmented systems of partial differential
equations in \cref{sec_charac_bifurc}. Then,
in \cref{sec_shape_opt}, we describe an algorithm for performing
PDE-constrained shape optimization constrained by the augmented system. In
\cref{sec_opt_bifurc_points}, we combine these numerical techniques in an
algorithm that modifies a shape to delay or advance a
given branch point. Finally, in \cref{sec_numer_ex}, we apply this
technique to delay the bifurcations of three numerical examples
involving the Allen--Cahn equation, the Navier--Stokes equations, and the
compressible hyperelasticity equation.

\section{Characterization of branch points} \label{sec_charac_bifurc}
Bifurcation diagrams of problems of the form $F(u,\lambda)=0$ may contain branch points (see~\cite{seydel2009practical}), which are
characterized by the change in the number of solutions as $\lambda$ varies in
a small interval around $\lambda^{\star}$. This implies that the Fr\'echet derivative $F_u(u^{\star},\lambda^{\star})$ 
with respect to the variable $u$ at the branch point $(u^{\star},\lambda^{\star})$ cannot be an isomorphism (otherwise,
by the implicit function theorem, there would exist locally a unique solution curve $(\lambda,u(\lambda))$,
but this is not compatible with \cref{def_bifurcation_points} below). Therefore, we are 
interested in the computation of singular points, i.e., solution
pairs $(u^{\star},\lambda^{\star})$ to \cref{eq:intro} such that $F_u(u^{\star},\lambda^{\star})$ does not have a bounded inverse.

The nonlinear operator $F$ may have a symmetry group $\mathcal{G}$ acting on $U$ such that for all group elements $g\in\mathcal{G}$ and $\lambda\in\mathbb{R}$, $F(u,\lambda)=0 \iff F(g\cdot u,\lambda)=0$. Here $g\cdot u$ represents the action of the element $g$ on the solution $u$. Such a group is called a Lie group when it has a structure of a smooth manifold and when the
group multiplication and inversion operations are smooth maps~\cite{golubitsky2012singularities,olver2000applications}. In the presence of such symmetries, one is typically interested
in points at which the number of \emph{distinct group orbits} varies, which are sets of the form: $\mathcal{G}(u)=\{g\cdot u \,: g\in\mathcal{G}\}$. In this work, we assume that no such Lie groups
exist for $F$, so that we may consider distinct solutions rather than distinct group orbits.

\begin{definition}[Branch point] \label{def_bifurcation_points}
Assume that the operator $F$ does not have a Lie group of symmetries. For $\lambda\in\mathbb{R}$, define
$\mathbf{U}_\lambda\coloneqq\{u\in U \,: F(u,\lambda)=0\}$
and, for $u^{\star}\in U$ and $\epsilon>0$, let
$B_\epsilon(u^{\star}) \coloneqq\{u\in U\,: \| u^{\star}- u\| <\epsilon\}$. A
branch point (also called bifurcation point) with respect to $\lambda$ is a
solution pair $(u^{\star},\lambda^{\star})$ to $F(u,\lambda)=0$ such that, for any
$\varepsilon>0$, the cardinality of $\mathbf{U}_\lambda\cap B_\varepsilon(u^{\star})$
changes when $\lambda$ passes $\lambda^{\star}$.
\end{definition}

\begin{figure}[htbp]
\vspace{0.2cm}
\centering
\begin{overpic}[width=0.3\textwidth]{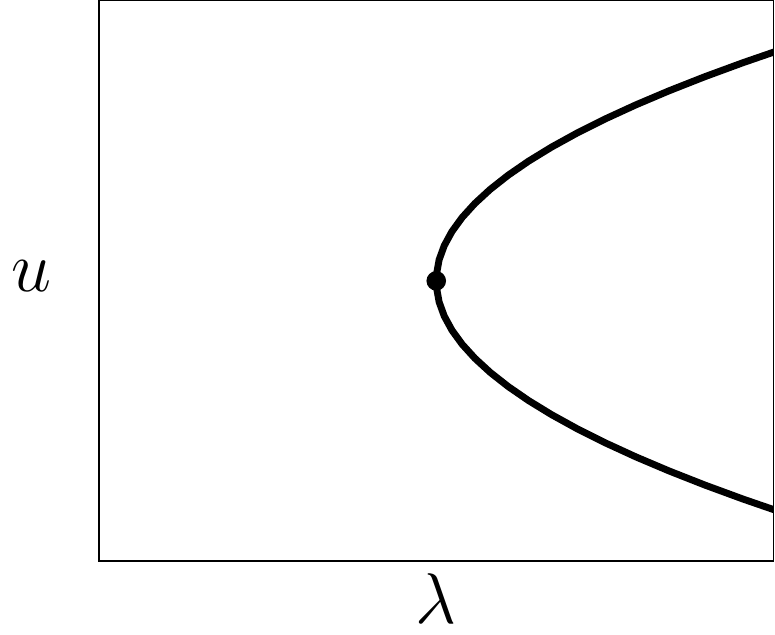}
\put(-2,80){(a)}
\end{overpic}
\hspace{0.3cm}
\begin{overpic}[width=0.3\textwidth]{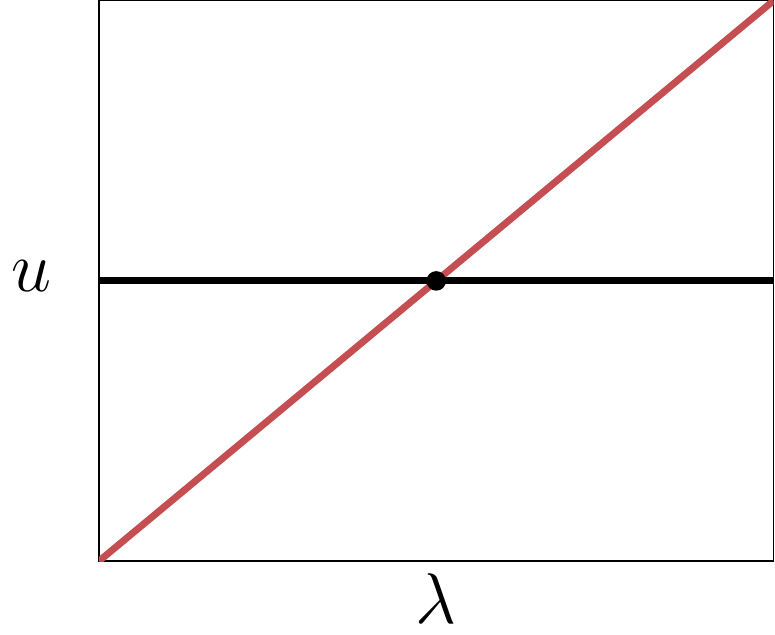}
\put(-2,80){(b)}
\end{overpic}
\hspace{0.3cm}
\begin{overpic}[width=0.3\textwidth]{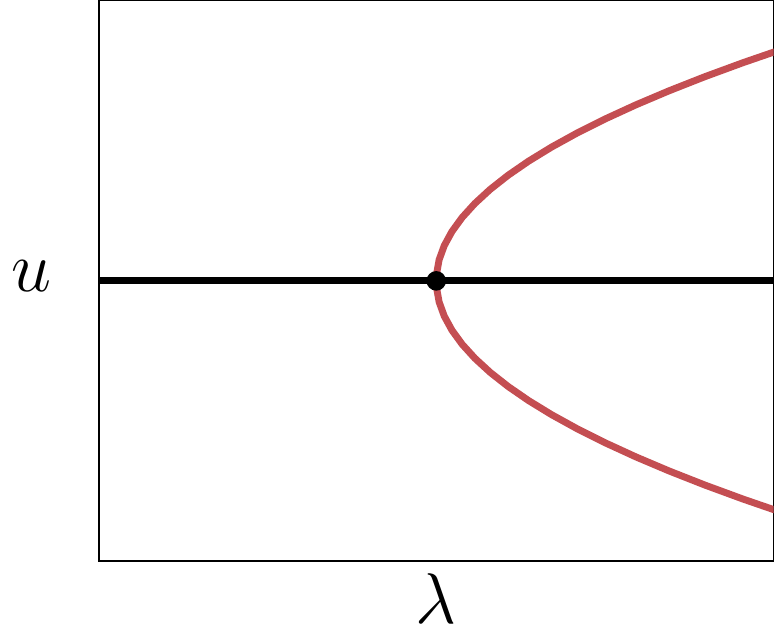}
\put(-2,80){(c)}
\end{overpic}
\caption{(a) Fold, (b) transcritical and (c) pitchfork bifurcations.}
\label{fig_bifurcation}
\end{figure}

We illustrate some different situations with branch points arising in the bifurcation diagram of a nonlinear PDE in \cref{fig_bifurcation}. \cref{fig_bifurcation}(a) features a fold bifurcation with a turning point. Panels (b) and (c) illustrate two potential scenarios for the intersection of two branches.
In a transcritical bifurcation (panel (b)), the trivial branch (indicated by a black line) may lose stability when one of the state eigenvalues crosses zero at the branch point, which results in an exchange of stability to the other branch, indicated by a red line. A pitchfork bifurcation (c) occurs when one solution gives birth to three, and often also involves an exchange of stability.

The idea of Seydel and Moore \& Spence for finding singular solutions to \cref{eq:intro} is to solve the following extended system~\cref{eq_moore_spence}:
find $(u,\lambda,\phi)\in U\times\mathbb{R}\times U$ such that
\begin{subequations} \label{eq_moore_spence}
\begin{align}
F(u,\lambda) &=0,\\
F_u(u,\lambda)\phi &= 0, \label{eq_MS_phi0} \\
\ell(\phi)-1 &=0, \label{eq_MS_phi} 
\end{align}
\end{subequations}
where $\ell:U\to\mathbb{R}$ is a functional chosen to
normalize $\phi$. Equations \cref{eq_MS_phi0} and \cref{eq_MS_phi} imply that $\phi$
is a non-zero eigenfunction associated to the zero eigenvalue of $F_u(u,\lambda)$.

\begin{remark}
Several choices of normalization conditions $\ell$ in the Moore--Spence system have been considered~\cite{cliffe1986use,winters1988onset}. 
In this work we set
\begin{equation*}
\ell(\phi) = \|\phi\|^2_U.
\end{equation*}
\end{remark}

\begin{remark}
In this work, we only consider simple singular points, those for which $\mathrm{dim}(\mathrm{ker}(F_u(u^{\star},\lambda^{\star}))) = 1$.
This includes nondegenerate turning points and simple symmetry-breaking pitchfork bifurcations~\cite{moore1980calculation,werner1984computation}. We further assume that the solutions
of \cref{eq_moore_spence} are isolated.
Bifurcation points of higher multiplicity, satisfying $\mathrm{dim}(\mathrm{ker}(F_u(u^{\star},\lambda^{\star}))) > 1$, also satisfy \cref{eq_moore_spence} \emph{a fortiori},
but may not be isolated.
It is possible to design augmented systems that are selective, in the sense that they admit certain kinds of branch points as solutions but not others
(e.g.~\cite{moore1980numerical,weber1981numerical}). These typically have the advantage that the solutions are isolated, but the augmented systems become larger and more complicated.
\end{remark}

\begin{remark}
Different generalizations have also been studied to compute Hopf bifurcations~\cite{strogatz2018nonlinear}, i.e.~points at which a periodic solution arises, in contrast to the stationary solutions considered in this work~\cite{griewank1983calculation,jackson1987finite,jepson1981numerical,roose1985direct}.
\end{remark}

Since we are interested in changing the shape of the domain $\Omega$ to control bifurcation diagrams, we now present a shape optimization technique constrained by partial differential equations.

\section{PDE-constrained shape optimization} \label{sec_shape_opt}

In this section, we give a brief introduction to PDE-constrained shape optimization.
A more detailed exposition of the mathematical and implementation techniques
employed here is given in \cite{paganini2018higher,paganini2020fireshape}.

A shape optimization problem is an optimization problem of the following
form: find a domain $\Omega\in \mathcal{U}$ that minimizes a shape functional
$J:\mathcal{U}\to\mathbb{R}$, where $\mathcal{U}$ denotes a set of admissible
domains. If evaluating the shape functional $J$ on a domain $\Omega$ requires
solving a PDE stated on $\Omega$, the shape optimization problem is said to
be PDE-constrained.

There are competing approaches to define the set of admissible domains
$\mathcal{U}$, such as the phase-field method \cite{BlGaFaHaSt14,GaHeHiKa15},
the level-set method \cite{AlJoTo02,La18} and the moving mesh method
\cite{AlPa06,paganini2018higher}, among others. In this work, we consider
the moving mesh method, because its neat integration with standard finite
element software allows the automation of several key steps in shape optimization
\cite{paganini2020fireshape}. The main idea of the moving mesh method is
to construct the set of admissible domains $\mathcal{U}$ by applying
diffeomorphisms to an initial domain $\Omega^0\subset\mathbb{R}^d$, that is,
\begin{equation*}
\mathcal{U} \coloneqq \{T(\Omega^0)\mid T\in\mathcal{T}\},
\end{equation*}
where $\mathcal{T}$ is the group of bi-Lipschitz isomorphisms on
$\mathbb{R}^d$, that is $T\in\mathcal{T}$ and $T^{-1}$ are Lipschitz continuous. 
In this setting, shape optimization is carried out by constructing a
geometric transformation $T$ so that $T(\Omega^0)$ minimizes $J$. This
approach is called moving mesh method because, in its lowest order
discretization, this method is equivalent to replacing the initial domain
$\Omega^0$ with a finite element mesh and optimizing the coordinates of the
mesh nodes while retaining the mesh connectivity \cite{paganini2018higher}.

\begin{remark}
The set of admissible domains contains also domains that are not Lipschitz.
If the shape functional $J$ is constrained to a PDE, then it is common to
restrict the domain of definition of $J$ to all elements of $\mathcal{U}$ that
are also Lipschitz to ensure that the PDE constraint is well defined. This
restriction can be obtained by further assuming, for instance, that
the transformations $T$ are continuously differentiable.
\end{remark}

Proving that a shape optimization problem admits a solution is not
straightforward. The main results rely on compactness arguments and require
additional assumptions on the collection of admissible domains $\mathcal{U}$
\cite[Ch. 6.2]{Al07}. Embedding these assumptions into numerical algorithms is
extremely challenging, to the point that in practice one contents oneself with
finding a sequence of domains along which the shape functional $J$ decreases.
In the framework of moving mesh methods, the sequence of domains is modelled via
a sequence of transformations that are commonly constructed iteratively. First,
one initializes $T^{(0)}$ with the identity transformation, that is, $T^{(0)}(x)
= I(x) \coloneqq x$ for every $x\in\mathbb{R}^d$. Then one constructs the
subsequent iterates using the composition formula
\begin{equation}\label{eq:T_update}
T^{(k+1)}=(I+\delta T^{(k)})\circ T^{(k)}\,, \quad k\geq 0,
\end{equation}
where $\delta T^{(k)}$ is a suitable geometry update. This sequence of iterates
corresponds to the sequence of domains $\{T^{(k)}(\Omega^0)\}$. More specifically,
the $(k+1)$-th element $T^{(k)}(\Omega^0)$ of this sequence of domains is obtained
by considering the previous domain $T^{(k)}(\Omega^0)$ and by adding to each of its
points $x\in T^{(k)}(\Omega^0)$ the corresponding vector $\delta T^{(k)}(x)$.

The geometry updates $\delta T^{(k)}$ can be computed in a standard fashion
using derivatives of the function $J$. As an example, in the next paragraph we
present a simple steepest descent shape optimization algorithm. More advanced
optimization strategies can be adapted to shape optimization in a similar way.

\begin{figure}[htbp]
\centering
\vspace{0.2cm}
\begin{overpic}[width=0.9\textwidth]{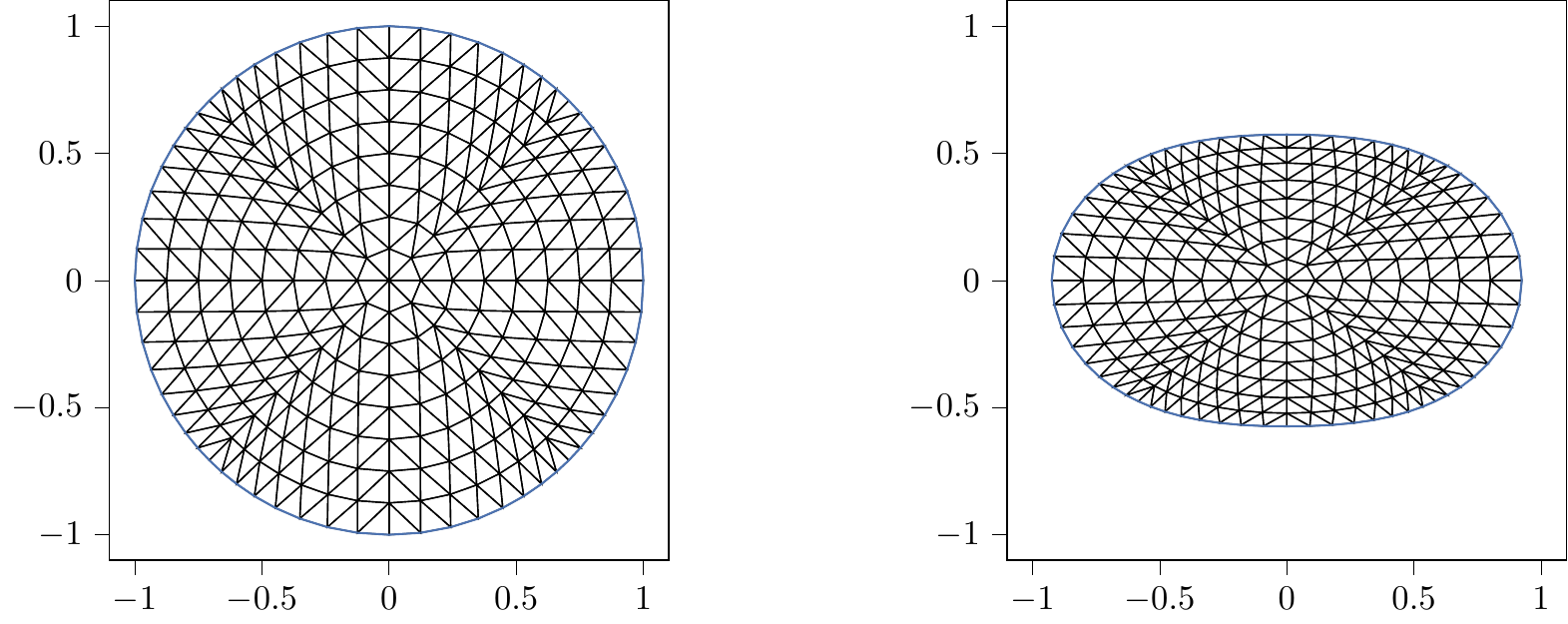}
\put(-2,38){(a)}
\put(55,38){(b)}
\end{overpic}
\caption{Initial domain (a) together with the modified domain after one
steepest descent update (b) computed using \cref{code1}. The optimal solution is an ellipse with major
and minor axes equal to 1 and 2/3, respectively. The objective functional
$J$ decreases from -0.59 to -1.01.}
\label{fig:minex}
\end{figure}

To construct the sequence $\{T^{(k)}\}$ using a steepest descent algorithm,
we compute the update $\delta T^{(k)}$ as the Riesz representative of the negative shape
derivative $-dJ$ evaluated at $T^{(k)}$ \cite{PaHi16,paganini2018higher}.
This means that $\delta T^{(k)}$ satisfies the equation
\begin{equation} \label{eq_Ries_rep}
(\delta T^{(k)},V)_\mathcal{H}=-dJ(T^{(k)}(\Omega^0);V) \quad \text{for all }
V\in \mathcal{H},
\end{equation}
where $\mathcal{H}$ is a suitable Hilbert space with inner product
$(\cdot,\cdot)_\mathcal{H}$, and 
\begin{equation} \label{eq_trans_diff}
dJ(T(\Omega^0);V) \coloneqq
\lim_{s\to 0^+}\frac{J((I+sV)\circ
T(\Omega^0))-J(T(\Omega^0))}{s}.
\end{equation}
Once $\delta T^{(k)}$ has been computed, the new transformation $T^{(k+1)}$ is
obtained using formula \cref{eq:T_update}, possibly scaling $\delta T^{(k)}$
by a sufficiently small step length $s \in (0, 1)$ to ensure that
$T^{(k+1)}$ is bijective~\cite[Lemma~6.13]{Al07}. \Cref{code1} shows how
to implement (in Firedrake) one step of the steepest descent method to minimize
the shape functional 
\begin{equation*}
J(\Omega) = \int_\Omega f\,\mathrm{d}x
\end{equation*}
for the specific two-dimensional function $f(x_1,x_2) =x_1^2 + (3 x_2/2)^2$. The
shape derivative of $J$ is
\begin{equation*}
dJ(\Omega; V) = \int_\Omega \nabla f\cdot V + f \nabla\cdot V\,\mathrm{d}x\,,
\end{equation*}
and the optimal solution is the zero-level set of $f$, which in this case
is an ellipse with major and minor axes equal to $1$ and $2/3$, respectively.
To compute the Riesz representative, we use piecewise linear Lagrangian finite
elements and consider the inner product
\begin{equation*}
(V, W)_{\mathcal{H}} = \int_\Omega DV \colon DW + V\cdot W\,\mathrm{d}x\,,
\end{equation*}
where $DV \colon DW$ is the Frobenius inner product
of the Jacobians of $V$ and $W$. \Cref{fig:minex} shows the initial mesh
and the mesh obtained after one optimization step. A similar code for a more
challenging shape functionals is given in \cite{ham2019automated}.

\begin{lstlisting}[basicstyle=\ttfamily, keywordstyle=\bfseries,frame=single,caption=Firedrake implementation of a steepest descent update for a shape optimization problem with a disk domain.,label=code1]
from firedrake import *

# Initial mesh and initial transformation
mesh = UnitDiskMesh(3)
T = mesh.coordinates.vector()

# Shape functional and its value on initial mesh
x, y = SpatialCoordinate(mesh)
f = x**2 + (3*y/2)**2 - 1 
J = f*dx
print("Initial J = %1.2f" % assemble(J))

# Solve (3.2) and update domain using (3.1)
H = mesh.coordinates.function_space()
dT = Function(H)
V = TestFunction(H)
A = inner(grad(dT),grad(V))*dx + dot(dT,V)*dx
dJ = dot(grad(f),V)*dx + f*div(V)*dx
solve(A + dJ == 0, dT) 
T += 0.5*dT.vector()

# Evaluate shape functional
print("Updated J = %1.2f" % assemble(J))
\end{lstlisting}

\begin{remark}
There are competing requirements on the choice of the Hilbert space
$\mathcal{H}$ \cite{paganini2018higher}. On the one hand, it is desirable
that the Riesz representative is sufficiently regular to ensure that $T^{(k)}$
is bi-Lipschitz. This suggest using an $H^{d/2+1}(\mathbb{R}^d)$-type space,
which is embedded in $W^{1,\infty}(\mathbb{R}^d)$~\cite[Sect.~5.6.3]{Ev10}.
On the other hand, it is desirable that \cref{eq_Ries_rep} can be solved
efficiently using standard finite elements, which suggest using an
$H^1(\mathbb{R}^d)$-type space. In this work, we follow the latter option
and employ an inner product based on the linear elasticity equations.
This is a standard choice in shape optimization~\cite{blauth2020nonlinear,etling2020first,schulz2016efficient} and, although it is not
guaranteed that $\delta T^{(k)}\in W^{1,\infty}(\mathbb{R^d})$
at the continuous level, its finite element approximation (based on standard Lagrangian finite elements)
belongs to $H^1(\Omega)\cap W^{1,\infty}(\Omega)$, and is thus suitable to
define the transformation $T^{(k)}$.
\end{remark}

For the numerical simulations in \Cref{sec_numer_ex} we employ the shape optimization software Fireshape
\cite{paganini2020fireshape}, which is based on the moving mesh method.
This choice
offers numerous convenient features, including automated shape differentiation
in the Unified Form Language~\cite{ham2019automated}, automated adjoint
computation via {pyadjoint}~\cite{dokken2020automatic,farrell2013automated}, integration with the finite element software Firedrake
\cite{rathgeber2016firedrake} for efficient linear solvers (e.g.~block preconditioners
and geometric multigrid) from PETSc~\cite{balay2020petsc}, and
access to numerous optimization algorithms via the Rapid Optimization Library~\cite{ridzal2017rapid}, including augmented Lagrangian algorithms~\cite[Chapt.~17.3]{nocedal2006numerical} to incorporate additional equality (and inequality) constraints such as volume or surface area
constraints~\cite[Sect.~2 and 3]{paganini2020fireshape}.

\section{Optimization of isolated branch points} \label{sec_opt_bifurc_points}

In light of \cref{sec_charac_bifurc,sec_shape_opt}, we formulate the problem of finding a shape of the
domain $\Omega^{\star}$ such that the operator $F$ admits a branch point
$(u^{\star},\lambda^{\star})$ at a target bifurcation parameter $\lambda^{\star}\in\mathbb{R}$
as follows:

\begin{equation} \label{eq_opt_moore_spence}
\min_{\Omega\in\mathcal{U}}
J(\Omega) \coloneqq (\lambda-\lambda^{\star})^2
\quad\textrm{subject to} \quad
\left\{\begin{aligned}
& F(u,\lambda)&=0\,,\\
& F_u(u,\lambda)\phi&=0\,,\\
& \Vert\phi\Vert_U^2 - 1&=0\,.
\end{aligned}\right.
\end{equation}
In this section, we describe an algorithm (\cref{alg_branch}) for solving \cref{eq_opt_moore_spence} and controlling branch points of nonlinear problems of the form of \cref{eq:intro}.

\renewcommand{\algorithmicrequire}{\textbf{Input:}}
\renewcommand{\algorithmicensure}{\textbf{Output:}}
\begin{algorithm}
\caption{Optimization of branch points}\label{alg_branch}
\begin{algorithmic}[1]
\Require Initial domain $\Omega^0$, initial guess $(\tilde{u},\tilde{\lambda})$ for the branch point  $(u^0,\lambda^0)$, target bifurcation parameter $\lambda^\star$, tolerance $\epsilon$
\Ensure Optimized domain $\Omega^{\star}$
\State Solve the eigenvalue problem~\cref{eq_eig_problem} to generate $n$ initial guesses $(\tilde{u},\tilde{\lambda},\phi_i)$
\For{$i=1:n$}
\State Solve the Moore--Spence system~\cref{eq_moore_spence} to obtain a solution $(\bar u_i, \bar
\lambda_i, \bar \phi_i)$
\EndFor
\State Select the initial solution $(u^0,\lambda^0,\phi^0)$ minimizing $\|\tilde{u}-\bar{u}_i\|_U$
\State $\lambda\gets \lambda^0$, $k\gets 1$
\While {$(\lambda-\lambda^{\star})^2>\epsilon$}
\State Evaluate the shape functional $J(\Omega)$ and compute shape update~\cref{eq_Ries_rep}
\State Update the mesh with~\cref{eq:T_update}
\State Solve the Moore--Spence system~\cref{eq_moore_spence} and find solution $(u^{(k)}, \lambda^{(k)}, \phi^{(k)})$
\If{the mesh is not tangled and~\cref{eq_close_sol} is satisfied}
\State $\lambda\gets \lambda^{(k)}$, $k\gets k+1$
\Else
\State Reject optimization step and decrease the optimization step size
\EndIf
\EndWhile
\end{algorithmic}
\end{algorithm}

Compared to classical PDE-constrained shape optimization, the problem
\cref{eq_opt_moore_spence} presents an additional difficulty: the state
constraint admits multiple solutions by design. More precisely, unless the
bifurcation diagram of $F(u, \lambda) = 0$ is trivial, there are values of
$\lambda$ for which the problem to find $u\in U$ such that $F(u,\lambda)=0$
admits multiple solutions $u$. Additionally, the bifurcation diagram of $F(u,
\lambda) = 0$ typically includes multiple branch points $(u,\lambda)$ (in
the sense of \cref{def_bifurcation_points}). 
For these reasons, to solve
\cref{eq_opt_moore_spence}, we must select not only an initial domain
$\Omega^0$, but also an initial isolated branch point $(u^0, \lambda^0)$ that we
wish to control. This selection can be performed by first computing the
bifurcation diagram of the initial domain $\Omega^0$ and selecting
$(u^0,\lambda^0)$ via visual inspection. Since most algorithms for
bifurcation analysis compute branch points approximately, we then refine the
approximation of $(u^0, \lambda^0)$ using the following procedure.  We select a
pair $(\tilde{u}, \tilde\lambda)\in U\times \mathbb{R}$ from the computed bifurcation diagram that both satisfies
$F(\tilde{u},\tilde\lambda)=0$ and that is sufficiently close to $(u^0,
\lambda^0)$. Then, we solve the following eigenvalue problem:
\begin{equation} \label{eq_eig_problem}
\text{find } (\mu,\phi) \text{ such that } F_{u}(\tilde u,\tilde \lambda)\phi=\mu \phi
\end{equation}
to determine the first $n\in\mathbb{N}$ (a typical value used is $n=5$) smallest
(single) eigenvalues in magnitude and corresponding eigenfunctions
$\phi_1,\ldots,\phi_n$ (normalized to $1$ with respect to the norm
$\|\cdot\|_U$).  The eigenvalue problem \cref{eq_eig_problem} is solved using a
Krylov--Schur algorithm with a shift-and-invert spectral
transformation~\cite{stewart2002krylov}, implemented in the SLEPc
library~\cite{hernandez2005slepc}. 
Finally, we form $n$ initial guesses $(\tilde u ,\tilde \lambda, \phi_i)$ and
solve the Moore--Spence system with Newton's method.
This generates a set of (at most) $n$ distinct solutions $(\bar u_i, \bar
\lambda_i, \bar \phi_i)$, from which we finally select the tuple minimizing the
difference between $\bar u_i$ and $\tilde{u}$ in the $U$-norm as the starting
point for shape optimization, that is, as $(u^0, \lambda^0, \phi^0)$.

\begin{figure}[htbp]
\vspace{0.2cm}
\centering
\begin{overpic}[width=0.5\textwidth]{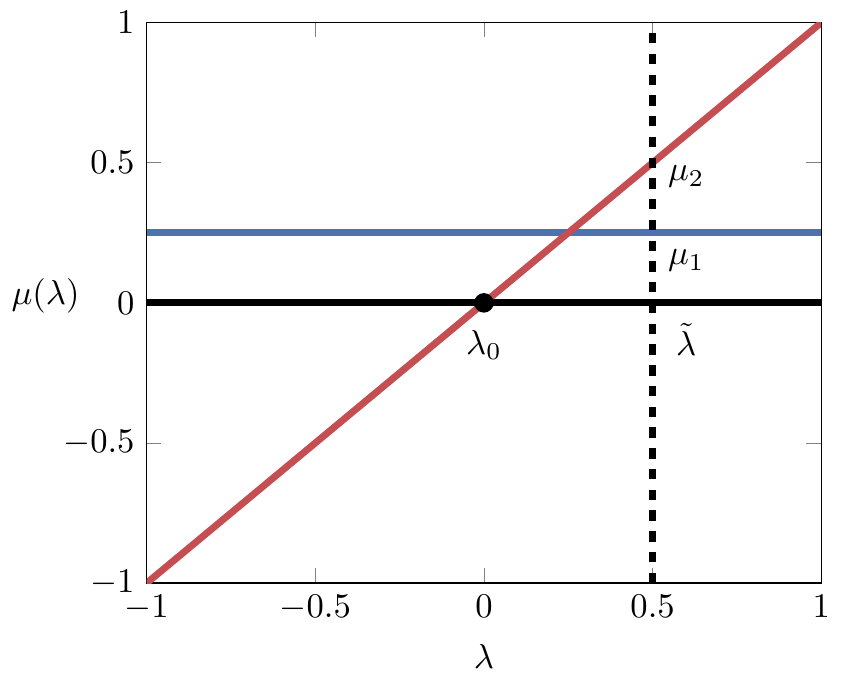}
\end{overpic}
\caption{Schematic of a situation where the eigenfunction associated with the second largest eigenvalue is a better initial guess that the smallest one for the Moore--Spence system. Two eigenvalue paths $\mu_1(\lambda)$ (blue) and $\mu_2(\lambda)$ (red) are plotted as a function of the parameter $\lambda$. At the solution of the Moore--Spence system $\lambda = \lambda_0$, $0 = \mu_2 < \mu_1$, but at the initial guess employed $\lambda = \tilde{\lambda}$, $\mu_2 > \mu_1$. One would therefore expect the \emph{second} lowest eigenfunction at $\lambda = \tilde{\lambda}$ to be a better initial guess.}
\label{fig_init}
\end{figure}

\begin{remark}
Using this procedure ensures that the initial state $(u^0,\lambda^0,\phi^0)$
is feasible. In practice, selecting $n>1$ is important, because situations like 
the one illustrated in \cref{fig_init} might arise; the eigenfunction associated with
the lowest eigenvalue for $\tilde{\lambda} \neq \lambda^0$ may not be a good
guess for the lowest eigenfunction at $\lambda^0$, since the eigenvalues may
swap order.
\end{remark}

Once $\Omega^0$ and $(u^0, \lambda^0, \phi^0)$ have been determined, we use
Fireshape \cite{paganini2020fireshape} to solve \cref{eq_opt_moore_spence}
iteratively with the trust-region optimization algorithm~\cite{conn2000trust}
and inner solver L-BFGS~\cite{liu1989limited} implemented in the ROL
optimization library~\cite{ridzal2017rapid}. To ensure that the computational
mesh is not tangled, we reject optimization steps for which the minimum of
the determinant of the Jacobian of the transformation $T$ defined in
\cref{eq:T_update} becomes negative. Additionally, we need to ensure that,
after a domain update, Newton's method converges to a solution $u^{(k+1)}$ to the
Moore--Spence system that lives on the same branch as the one computed at the
previous iteration, i.e.~is close enough to the solution $u^{(k)}$ obtained 
in the previous iteration. Hence, undesired branch-switching, caused by multiple solutions to the state equation, might occur during the shape optimization algorithm as observed in~\cite{williams2021shape}. To handle this, we reject domain updates that do not satisfy the inequality 
\begin{equation} \label{eq_close_sol}
\| u^{(k)} \circ (I+\delta T^{(k)})^{-1} - u^{(k+1)}\|_U \leq C \| u^{(k+1)} \|_U,
\end{equation}
where $I+\delta T^{(k)}$ is the diffeomorphism that updates the domain geometry in the
$k$-th iteration, and $C$ is an arbitrary constant ($C=0.1$ in
\cref{sec_allen,sec_beam}, $C=1$ in \cref{sec_ns}).

\section{Numerical examples} \label{sec_numer_ex}

In this section, we apply the shape optimization strategy
described in \cref{sec_opt_bifurc_points} to three test cases based on three
different differential equations $F$: the Allen--Cahn equation,
the Navier--Stokes equations, and the hyperelasticity equations. The bifurcation
diagrams displayed in this section are computed with deflated
continuation~\cite{farrell2016computation,farrell2015deflation}, complemented
with arclength continuation. These
techniques have been successfully applied to a wide range of physical
problems such as the deformation of a hyperelastic
beam~\cite{farrell2018a}, liquid
crystals~\cite{emerson2018computing,xia2021a}, Bose--Einstein condensates~\cite{boulle2020deflation,charalampidis2020bifurcation,charalampidis2018computing},
and fluid dynamics~\cite{boulle2021bifurcation}. However, our optimization
strategy is not tied to a specific algorithm for bifurcation analysis, and
other options may be used.

\subsection{Allen--Cahn equation} \label{sec_allen}

As first test case, we consider the cubic-quintic
Allen--Cahn equation, which can be used to model phase
separation~\cite{allen1979microscopic,uecker2014pde2path}.
In strong form, the Allen--Cahn equation is given by
\begin{subequations} \label{eq_ac_dynamics_systems}
\begin{align}
F(u,\lambda)\coloneqq-0.25\nabla^2u -\lambda u - u^3 +u^5 &= 0 \quad \text{in }\Omega,\\
u&=0 \quad \text{on }\partial\Omega.
\end{align}
\end{subequations}
In this example, we consider two-dimensional domains $\Omega \subset \mathbb{R}^2$. Therefore, we can set $U = H^1_0(\Omega)$. We use two different initial domains: a square $\Omega_\textrm{square}$
with rounded corners and edge length 2 (see
\cref{fig_allen_cahn_opt_square}(a)), and a unit disk $\Omega_\textrm{disk}$.
The computational meshes of
$\Omega_\textrm{square}$ and $\Omega_\textrm{disk}$ are both generated using Gmsh~\cite{geuzaine2009gmsh} and are composed of $9070$ and $13571$ triangles, respectively. The Allen--Cahn equation is discretized using piecewise linear continuous Lagrangian finite elements.

\begin{figure}[htbp]
\vspace{0.2cm}
\begin{center}
\begin{overpic}[width=\textwidth]{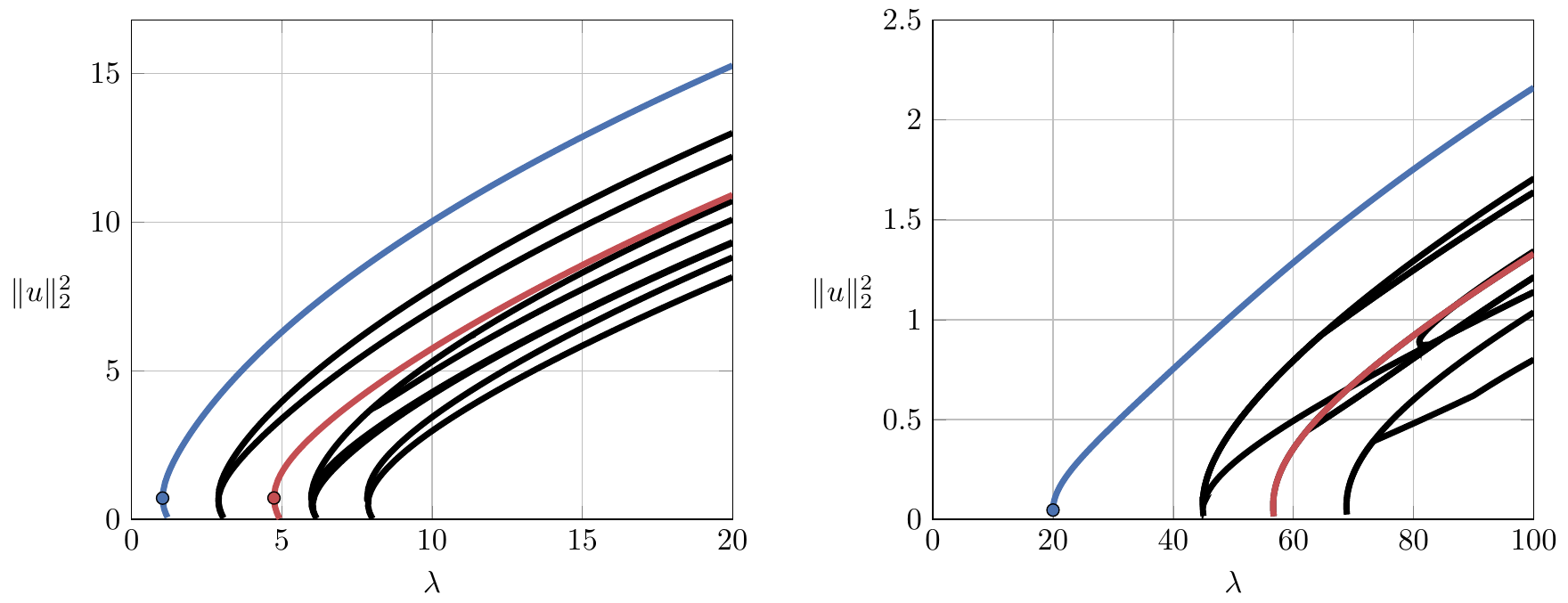}
\put(0,35.5){(a)}
\put(51,35.5){(b)}
\end{overpic}
\end{center}
\caption{Bifurcation diagrams of the Allen--Cahn problem defined on a square with round corners $\Omega_\textrm{square}$ (a), together with the diagram after delaying the first bifurcation to $\lambda=20$ (b). The first and third branches are respectively highlighted in blue and red.}
\label{fig_allen_cahn}
\end{figure}

\begin{figure}[htbp]
\vspace{0.2cm}
\begin{center}
\begin{overpic}[width=\textwidth]{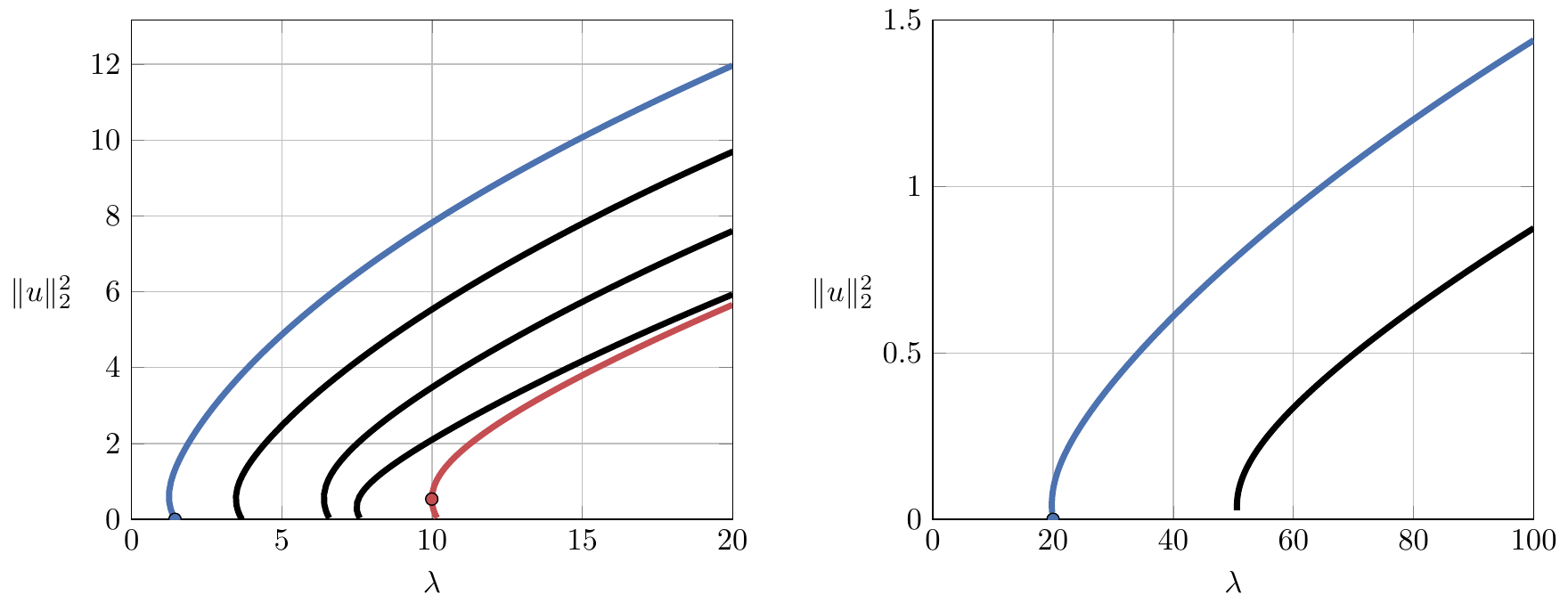}
\put(0,35.5){(a)}
\put(51,35.5){(b)}
\end{overpic}
\end{center}
\caption{Bifurcation diagrams of the Allen--Cahn problem defined on the unit disk $\Omega_\textrm{disk}$ (a), together with the diagram after delaying the first bifurcation to $\lambda=20$ (b). The first and fifth branches are respectively highlighted in blue and red.}
\label{fig_allen_cahn_disk}
\end{figure}

\begin{figure}[htbp]
\vspace{0.2cm}
\begin{center}
\begin{overpic}[width=\textwidth]{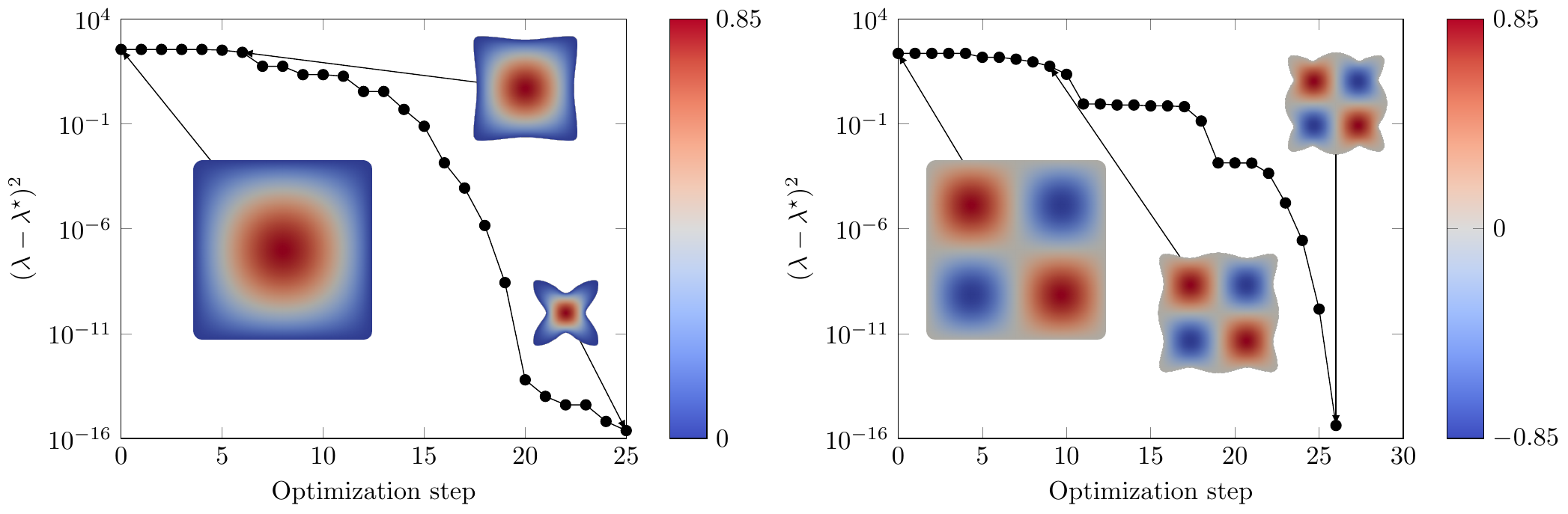}
\put(0,33.5){(a)}
\put(49.5,33.5){(b)}
\end{overpic}
\end{center}
\caption{Evolution of the domain during the shape optimization of the first (a) and third (b) branches (highlighted in blue and red in \cref{fig_allen_cahn}(a)) to the Allen--Cahn equation with a target branch point of $\lambda^{\star}=20$, initialized with the square domain $\Omega_\textrm{square}$. The curve represents the value of the objective function $(\lambda-\lambda^{\star})^2$ at each optimization step of the algorithm. The domains are plotted to the same scale.}
\label{fig_allen_cahn_opt_square}
\end{figure}

\begin{figure}[htbp]
\vspace{0.2cm}
\begin{center}
\begin{overpic}[width=\textwidth]{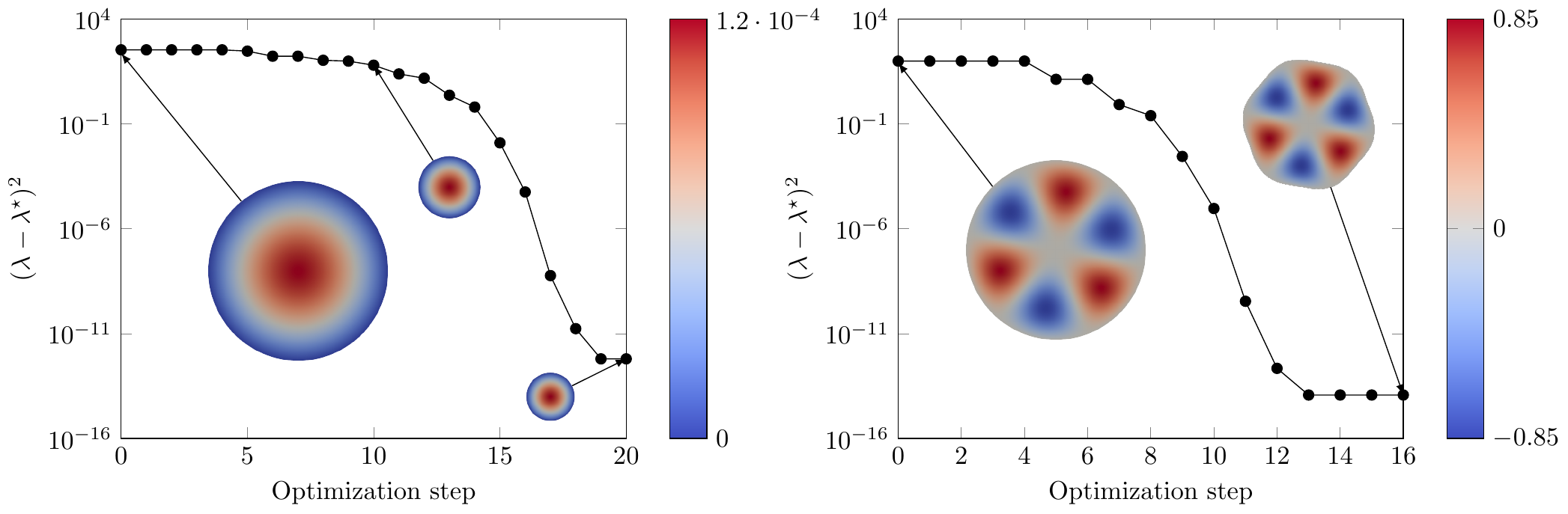}
\put(0,33.5){(a)}
\put(49.5,33.5){(b)}
\end{overpic}
\end{center}
\caption{Evolution of the domain during the shape optimization of the first (a) and fifth (b) branches (highlighted in blue and red in \cref{fig_allen_cahn_disk}(a)) to the Allen--Cahn equation with a target branch point of $\lambda^{\star}=20$, initialized with the disk domain $\Omega_\textrm{disk}$. The curve represents the value of the objective function $(\lambda-\lambda^{\star})^2$ at each optimization step of the algorithm. The domains are plotted to the same scale.}
\label{fig_allen_cahn_opt_disk}
\end{figure}

In \cref{eq_ac_dynamics_systems}, several bifurcations occur as $\lambda$
varies. In the left panels of \cref{fig_allen_cahn,fig_allen_cahn_disk}, we
show the bifurcation diagrams computed with deflated
continuation~\cite{farrell2016computation} for $\lambda\in[0,20]$ when the
computational domain is $\Omega_\textrm{square}$ and $\Omega_\textrm{disk}$,
respectively. We notice that several branches arise in both cases from the
ground state $u=0$. We can order these branches by numbering them
from the left. We observe that the solutions $u$ associated to higher
branches feature more complicated patterns than the states belonging to the
first branches (see
\cref{fig_allen_cahn_opt_square,fig_allen_cahn_opt_disk}). Another
interesting observation is the presence a fold bifurcation on each
branch, due to the cubic-quintic terms~\cite{uecker2014pde2path}.

We now apply the optimization method described in
\cref{sec_opt_bifurc_points} to delay the birth of the first branch from the
ground state in the Allen--Cahn equation to a target bifurcation parameter of
$\lambda^{\star}=20$. Here, the whole domain is deformable. In \cref{fig_allen_cahn_opt_square,fig_allen_cahn_opt_disk},
we display the evolution of the optimization functional
$\mathcal{R}(\lambda)=(\lambda-\lambda^{\star})^2$ with respect to the trust-region
iterations. We also show some snapshots of the solution $u$ and domain $T(\Omega^0)$ for certain
values of $\lambda$. We observe that the optimization algorithm reduces the
value of the objective functional $\mathcal{R}$ to less than $10^{-11}$
(\cref{fig_allen_cahn_opt_square,fig_allen_cahn_opt_disk}) in about 20
optimization steps.

To verify that the returned solutions satisfy the original objective,
we recompute the bifurcation diagrams (using deflated and arclength
continuation) on the optimized shapes. The 
right panels of \cref{fig_allen_cahn,fig_allen_cahn_disk} illustrate the
resulting bifurcation diagrams using the optimized shapes, starting from
$\Omega^0 = \Omega_{\textrm{square}}$ and $\Omega^0 = \Omega_{\textrm{disk}}$, respectively.
In both cases, we observe that the bifurcation structure is controlled
in the desired way, and that
the first branch arises at the desirable value of the bifurcation parameter.
As expected, the new shape also affects the birth of the later branches. For
instance, the second branch of the diagram with $\Omega_{\textrm{disk}}$
originally arises at $\lambda\approx 4$ and gets delayed to $\lambda\approx
50$ in \cref{fig_allen_cahn_disk}, while the higher branches are still not
visible in the range $\lambda\in[0,100]$. 

This example also showcases the
complexity and nonlinearity of the problem studied in this work. None of the
problems at any level have a unique solution: the PDE may have
multiple solutions for a given bifurcation parameter $\lambda$;
the Moore--Spence system has multiple solutions due to the presence of
several branch points; and the shape optimization problem admits
multiple solutions, as found from different initial shapes. The shapes displayed in
\cref{fig_allen_cahn_opt_square}(a) and \cref{fig_allen_cahn_opt_disk}(a) are
both solutions to the same shape optimization problem but initialized using
$\Omega_{\textrm{square}}$ and $\Omega_{\textrm{disk}}$.

Finally, we demonstrate that the same approach can be employed to control arbitrary
branches. We repeat the previous experiment and shape optimize
$\Omega_{\textrm{square}}$ to postpone the birth of the third branch (which
arises near $\lambda = 5$, see the left panel of \cref{fig_allen_cahn}) to
$\lambda^{\star}=20$. In \cref{fig_allen_cahn_opt_square}(b), we observe that the
shapes obtained by this procedure, the value of $\mathcal{R}$ falls below
machine precision in roughly 26 optimization steps. Similarly, in \cref{fig_allen_cahn_opt_disk} (b) we
show the optimization history when delaying the birth of the fifth
bifurcation of $\Omega_\textrm{disk}$ from (roughly) $\lambda=10$ to
$\lambda^{\star}=20$. In this case, the optimized shape resembles a (smoothened)
hexagon with edges at the extrema of the function $u$ (indicated by blue and
red colours).

\subsection{Navier--Stokes equations} \label{sec_ns}

In the second test case, we consider the functional $F$ associated to the
nondimensionalized steady incompressible Navier--Stokes equations. These read
\begin{align*}
-\nabla\cdot \left(\frac{2}{\textrm{Re}}\epsilon(u)\right)+u\cdot\nabla u+\nabla p&=0\quad \text{in }\Omega,\\
\nabla\cdot u&=0\quad \text{in }\Omega,\\
u&=0 \quad \text{on } \partial\Omega\backslash(\Gamma_{\text{in}}\cup \Gamma_{\text{out}}),\\
u&=g \quad \text{on } \Gamma_{\text{in}},\\
p\vec{n}-\frac{2}{\textrm{Re}}\epsilon(u)\cdot \vec{n} &= 0  \quad \text{on } \Gamma_{\text{out}},
\end{align*}
where $u$ is the fluid velocity, $\epsilon(u)=\frac{1}{2}(\nabla u+\nabla
u^\top)$, $p$ is the pressure, and Re is the Reynolds number, which in this
example plays the role of bifurcation parameter. In particular, we consider
an initial geometry formed by the union of two rectangles union of two
rectangles, $\Omega=[0,2.5]\times[-1,1]\cup [2.5,150]\times [-6,6]$, and
represents a channel subject to sudden expansion. The inflow is modelled with
a Poiseuille flow $g(x,y) = (1-y^2,0)^\top$ at the leftmost inlet boundary $\Gamma_\text{in}$, whereas
an outflow boundary condition is imposed on the rightmost outlet boundary  $\Gamma_\text{out}$. The
remaining parts of the boundary carry no-slip boundary conditions. The velocity and pressure are discretized with the inf-sup stable
Taylor--Hood finite element, i.e.~the vector-valued piecewise quadratic continuous Lagrangian element for velocity,
and piecewise linear continuous Lagrangian element for pressure.

\begin{figure}[htbp]
\vspace{0.2cm}
\begin{center}
\begin{overpic}[width=\textwidth]{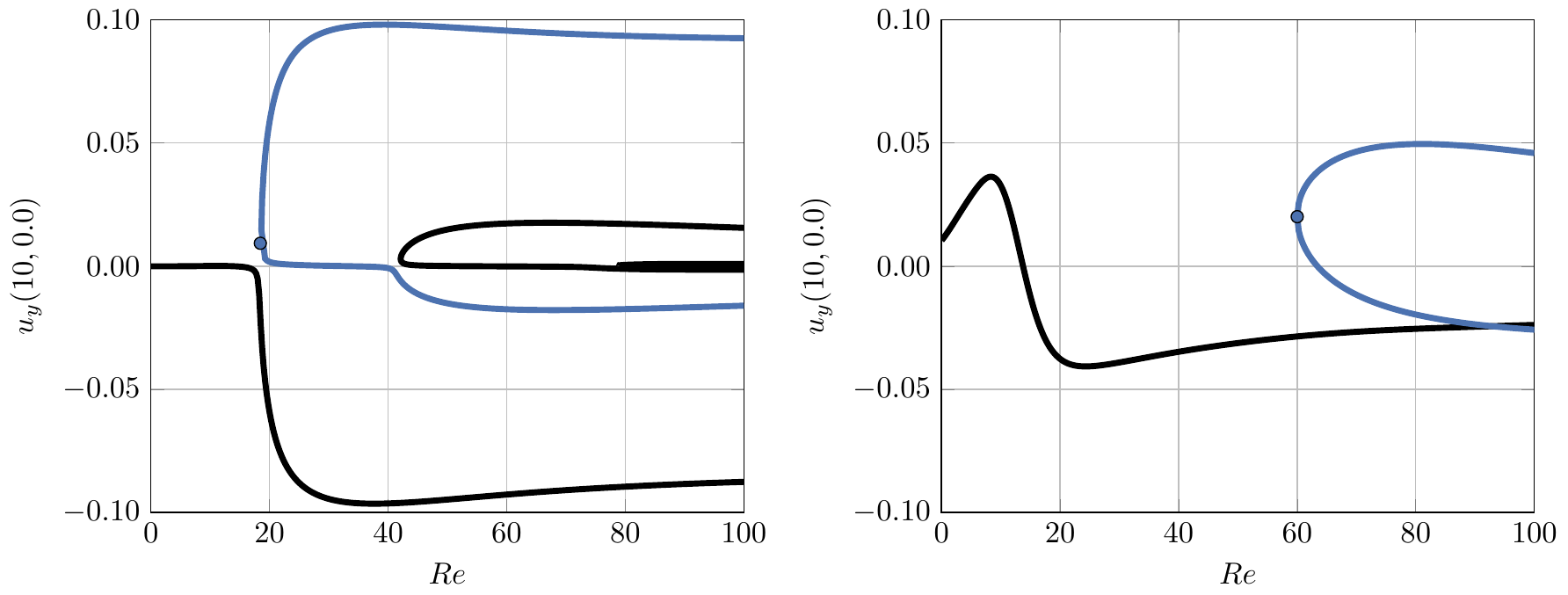}
\put(0,35.5){(a)}
\put(51,35.5){(b)}
\end{overpic}
\end{center}
\caption{Left: Initial bifurcation diagram of the Navier--Stokes equations defined on the channel domain $\Omega=[0,2.5]\times[-1,1]\cup [2.5,150]\times [-6,6]$. Right: Bifurcation diagram computed on the optimized channel displayed in \cref{fig_opt_navier_stokes}(c), where the first bifurcation branch (highlighted in blue) arises at $\textrm{\normalfont Re}=60$ through a fold bifurcation. The two diagrams are computed with deflated continuation and arclength continuation using $u_y(10,0)$ as a diagnostic, where $u_y$ is the vertical component of the velocity field $u$.}
\label{fig_navier_stokes}
\end{figure}

It is well known that this problem admits multiple solutions depending on the
values of the bifurcation parameter Re~\cite{fearn1990}.  The bifurcation
diagram for the initial mesh is shown in \cref{fig_navier_stokes}(a). The birth
of the first bifurcation from the symmetric ground state occurs at
$\textrm{Re}\approx 18.7$.  If the mesh were $\mathbb{Z}_2$-symmetric about $y =
0$, the bifurcation diagram would be symmetric and connected, with branches
originating at pitchfork bifurcations from the symmetric ground state. Since
symmetry is broken by the use of the unstructured mesh, the pitchforks are
perturbed into disconnected fold bifurcations. This phenomenon is very difficult
to avoid in practice (especially on complex industrial geometries), but
deflation deals with it robustly.

\begin{figure}[htbp]
\vspace{0.2cm}
\begin{center}
\begin{overpic}[width=0.9\textwidth]{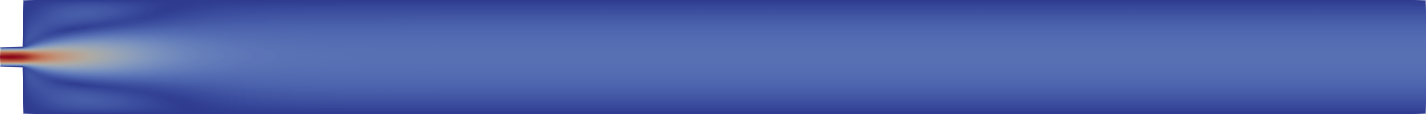}
\put(-5,3.3){(a)}
\end{overpic}\\
\vspace{0.4cm}
\begin{overpic}[width=0.9\textwidth]{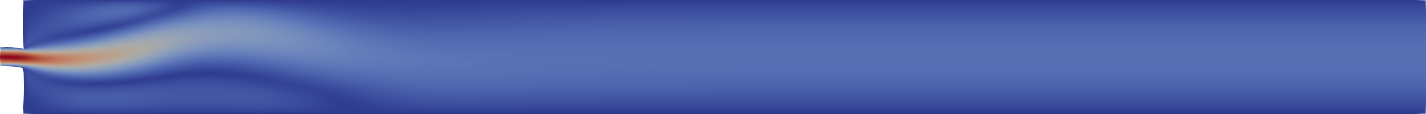}
\put(-5,3.3){(b)}
\end{overpic}\\
\vspace{0.4cm}
\begin{overpic}[width=0.9\textwidth]{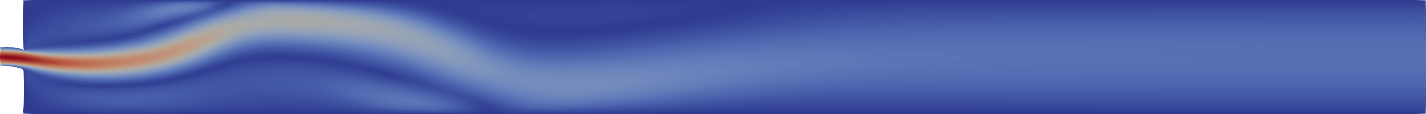}
\put(-5,3.3){(c)}
\end{overpic}\\
\vspace{0.4cm}
\begin{overpic}[width=0.9\textwidth]{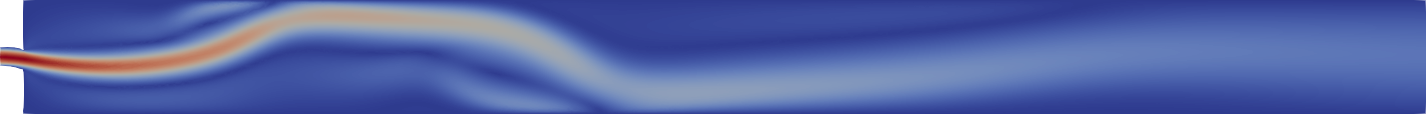}
\put(-5,3.3){(d)}
\end{overpic}\\
\vspace{0.4cm}
\begin{overpic}[width=0.4\textwidth]{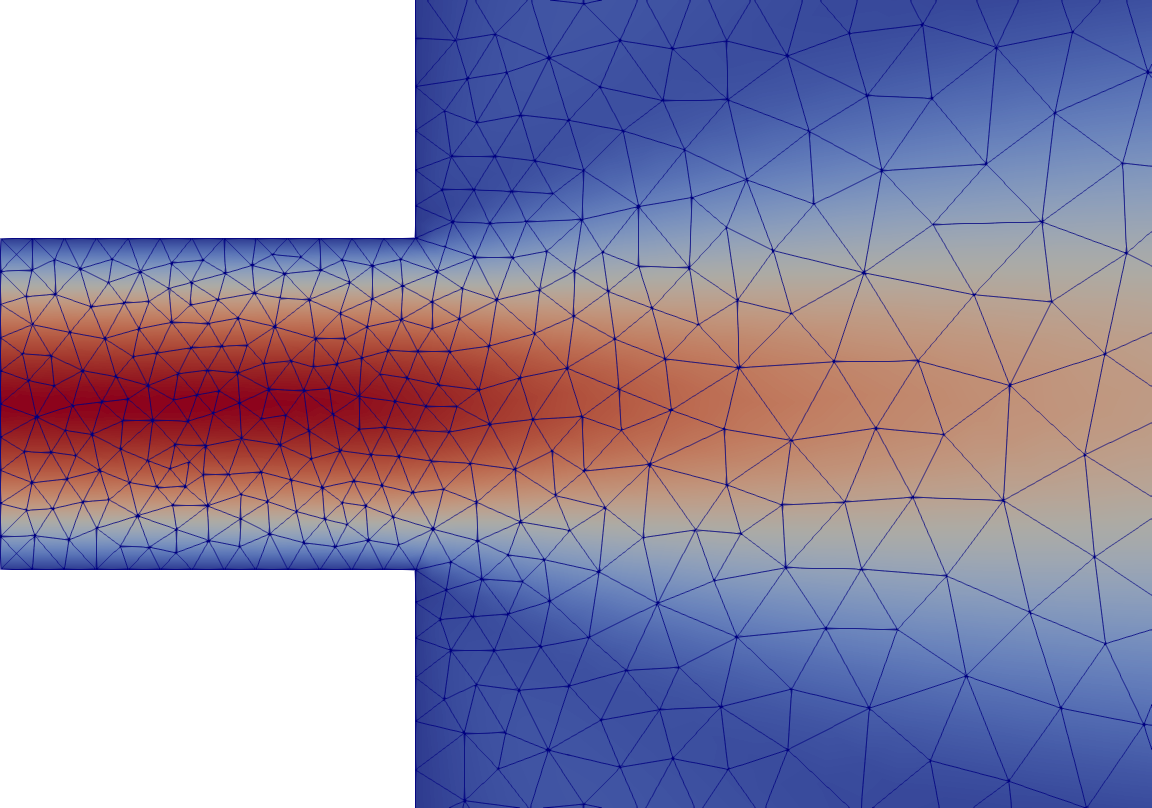}
\put(0,65){(a)}
\end{overpic}
\hspace{1.0cm}
\begin{overpic}[width=0.4\textwidth]{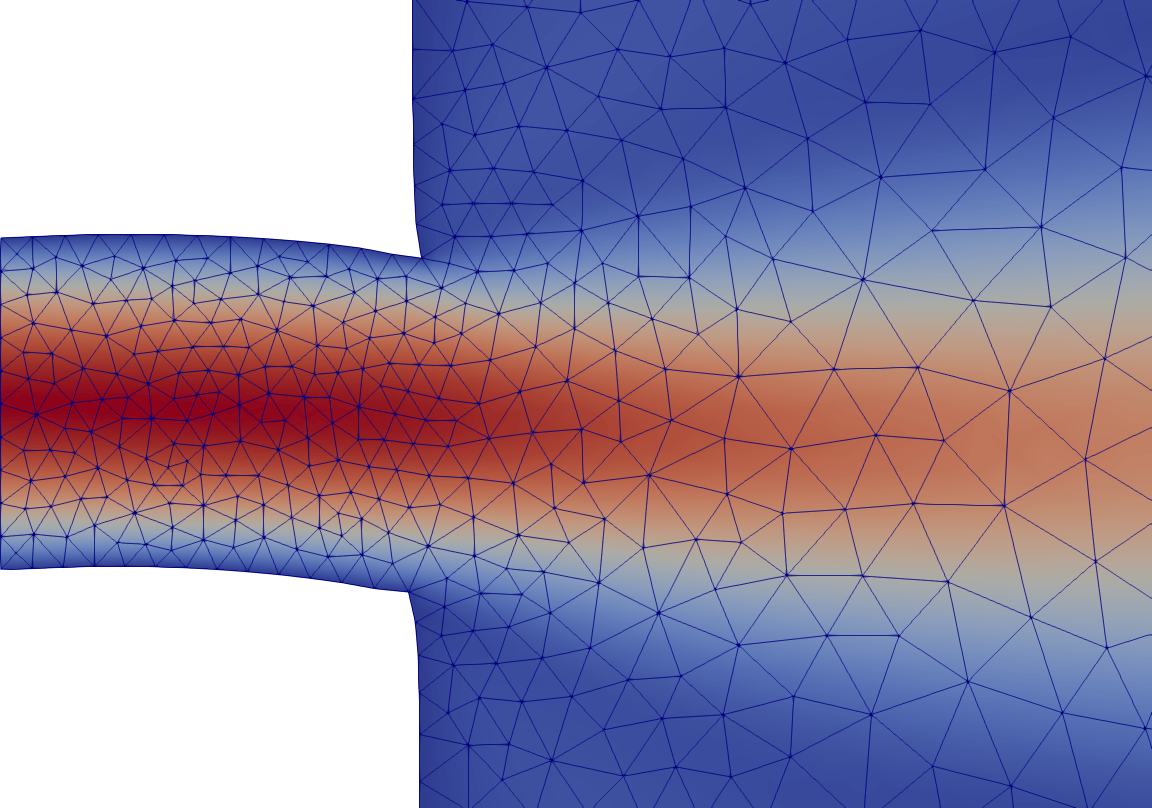}
\put(0,65){(b)}
\end{overpic}\\
\vspace{0.4cm}
\begin{overpic}[width=0.4\textwidth]{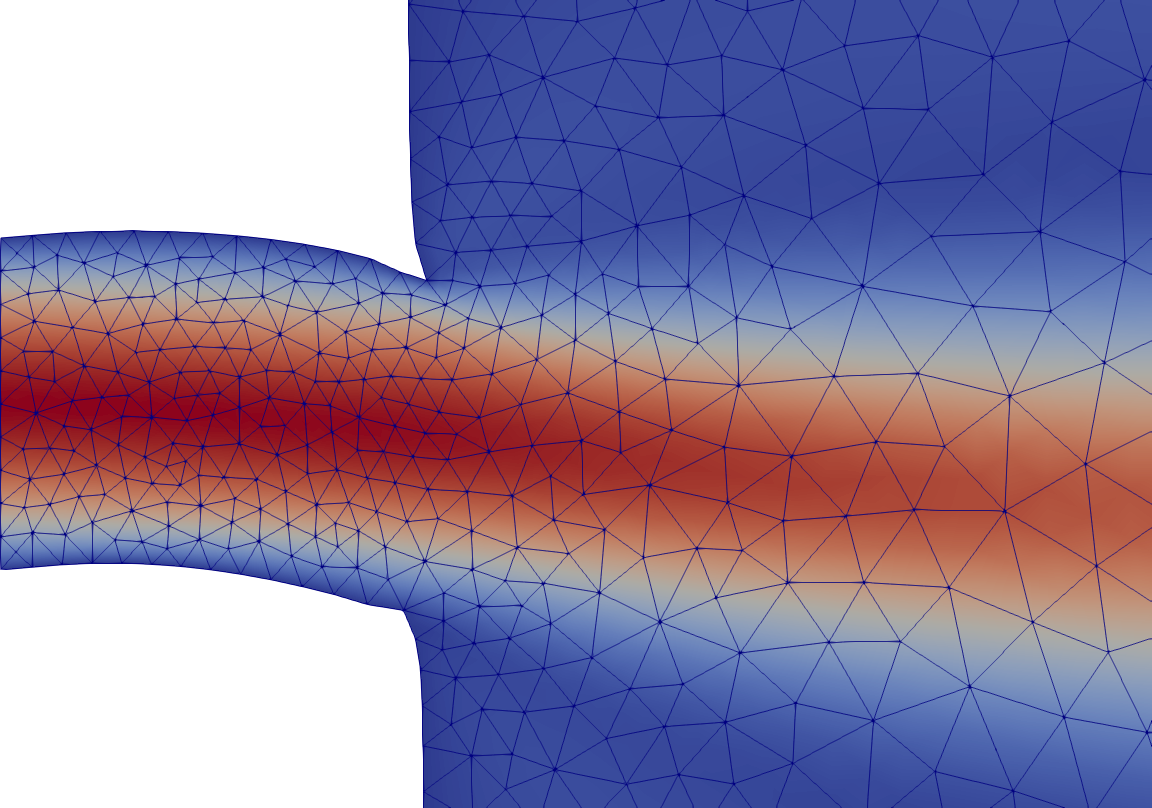}
\put(0,65){(c)}
\end{overpic}
\hspace{1.0cm}
\begin{overpic}[width=0.4\textwidth]{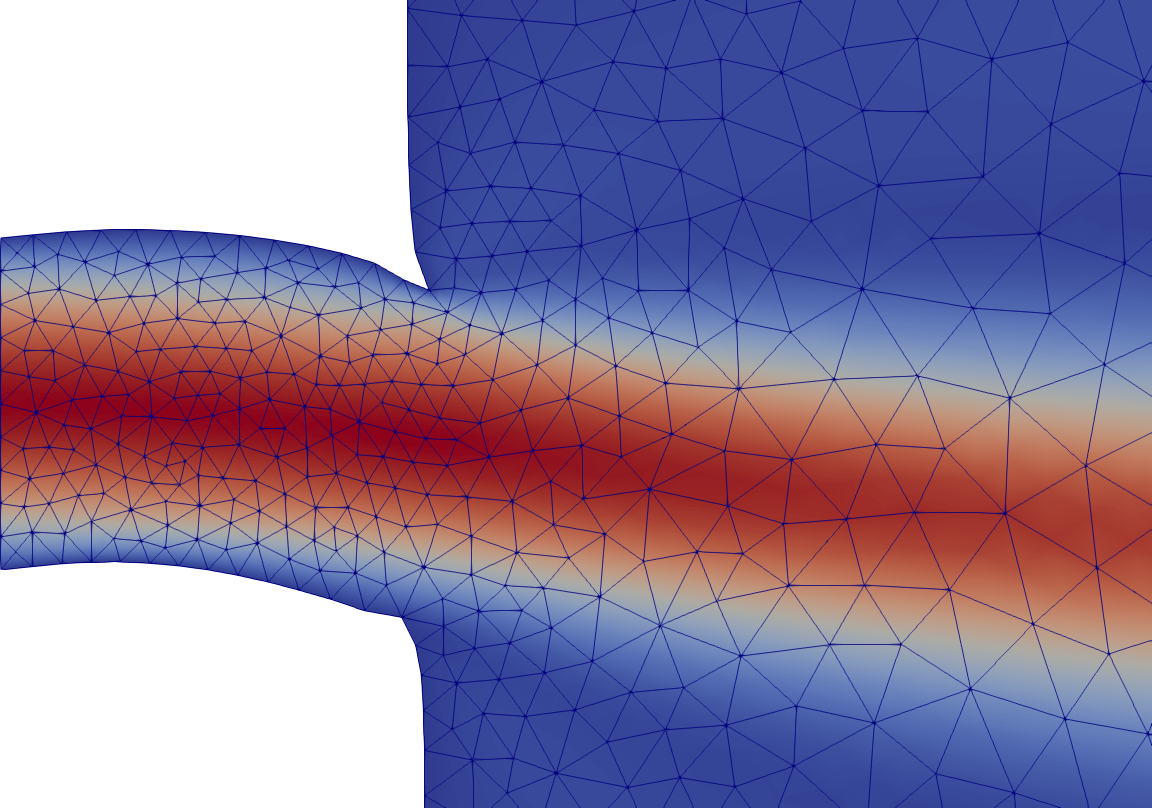}
\put(0,65){(d)}
\end{overpic}\\
\vspace{0.4cm}
\begin{overpic}[width=0.4\textwidth]{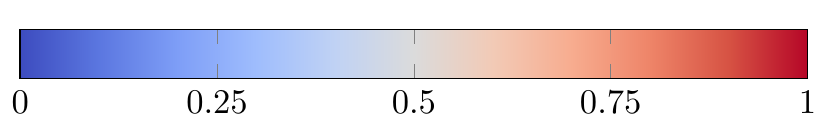}
\end{overpic}
\end{center}
\caption{Velocity magnitude of the solution to the Navier--Stokes equations at the first branch point arising at $\textrm{\normalfont Re} = $ 18.7 (a), 30 (b), 60 (c), 100 (d), together with a magnification of the transformed mesh at the inlet.}
\label{fig_opt_navier_stokes}
\end{figure}

Using the shape optimization approach described in
\cref{sec_opt_bifurc_points}, we modify the shape of the inlet to delay the
birth of the first bifurcation from the ground state solution.
Since we are interested in changing the shape of the inlet, 
we fix all boundaries during the shape optimization algorithm
except for the following segments:
\[[(0,1),(2.5,1)],\quad [(2.5,1),(2.5,6)],\quad [(0,-1),(2.5,-1)],\quad [(2.5,-1),(2.5,-6)].\]

\cref{fig_opt_navier_stokes}(a) displays the velocity
magnitude of the solution on the initial mesh at the bifurcation
point located at $Re\approx 18.7$. In \cref{fig_opt_navier_stokes}(b)-(d), we
report the shape of the inlet channel after delaying the branch point to
$Re=30,60,100$, respectively. We observe that the symmetry of the solution
illustrated in \cref{fig_opt_navier_stokes}(a) is broken as we delay
the birth of the bifurcation. The shape of the inlet channel
becomes bent and forces the fluid to flow downward as the value of the
branch point increases from $Re\approx 18.7$ to $Re = 100$.
In each case, the shape optimization algorithm required between 20 and 30 steps 
to minimize the loss function to nearly machine precision. 

Finally, we apply deflated continuation to compute the bifurcation diagram on
the bent domain shown in
\cref{fig_opt_navier_stokes}(c), for which the first branch point
happens at $Re=60$. The resulting bifurcation diagram is displayed in
\cref{fig_navier_stokes}(b). We observe that the first bifurcation branch
effectively arises at $Re=60$ through a fold bifurcation and is disconnected
from the original solution branch.

\subsection{Hyperelastic beam} \label{sec_beam}

We consider the shape optimization of a hyperelastic
beam under compression. The potential energy of the beam is given by
\begin{equation} \label{eq_energy_hyper}
\Pi(u) = \int_\Omega\psi(u)\,\mathrm{d}x - \int_\Omega B\cdot u\,\mathrm{d}x, 
\end{equation}
where $\Omega=[0,1]\times[0,0.1]$ is the initial domain,
$u:\Omega\to\mathbb{R}^2$ is the unknown displacement, $B=(0,-1000)^\top$ is a
constant body force per unit area, and $\psi$ is a Neo--Hookean strain energy density
given by
\begin{equation}\label{eq:strainenergydensity}
\psi=\frac{\mu}{2}(I_C-2)-\mu\ln(J)+\frac{\lambda}{2}\ln(J)^2.
\end{equation}
Here $\mu=E/(2(1+\nu))\approx 3.8\times 10^6$ and $\lambda=E\nu((1+\nu)(1-2\nu))\approx 5.8\times 10^6$
denote the Lam\'e parameters (chosen to yield
a Young's modulus $E = 10^6$ and Poisson ratio $\nu = 0.3$), $J=\textrm{det}(F)$, $F=I+\nabla u$,
$I_C=\textrm{tr}(C)$, and $C = F^\top F$.
Finally, we complement the potential energy \cref{eq_energy_hyper}
with the Dirichlet boundary conditions 
\[u(0,\cdot)=(0,0)^\top, \qquad u(1,\cdot)=(-\epsilon,0)^\top,\]
where the variable load $\epsilon\in[0,0.2]$ plays the role of bifurcation
parameter. On the remaining part of the boundary, denoted by $\Gamma_N$, we impose natural traction-free boundary conditions. 

Let $n$ denote the outward normal vector and
$P(u)$ be the first Piola--Kirchhoff stress tensor, calculated as the
derivative of the stored energy density function:
\[P(u)=\frac{\partial W}{\partial F} = \mu F(\tr(C)I)-\mu F^\top +
\frac{\lambda}{2}\ln(J)F^\top,\] where $W(F)\coloneqq\psi(u)$. 
The boundary value problem considered in this section is
\begin{alignat*}{3}
-\nabla\cdot P(u) &= B\quad &&\text{in }\Omega,\\
u&=0\quad &&\text{on }\Gamma_L,\\
u&=(-\epsilon,0)^\top\quad &&\text{on }\Gamma_R,\\
P(u)\cdot n&=0\quad &&\text{on }\Gamma_N,
\end{alignat*}
where $\Gamma_L \coloneqq \{0\}\times [0,0.1]$ and $\Gamma_R\coloneqq \{1\}\times [0,0.1]$ denote the left and right boundaries of $\Omega$.

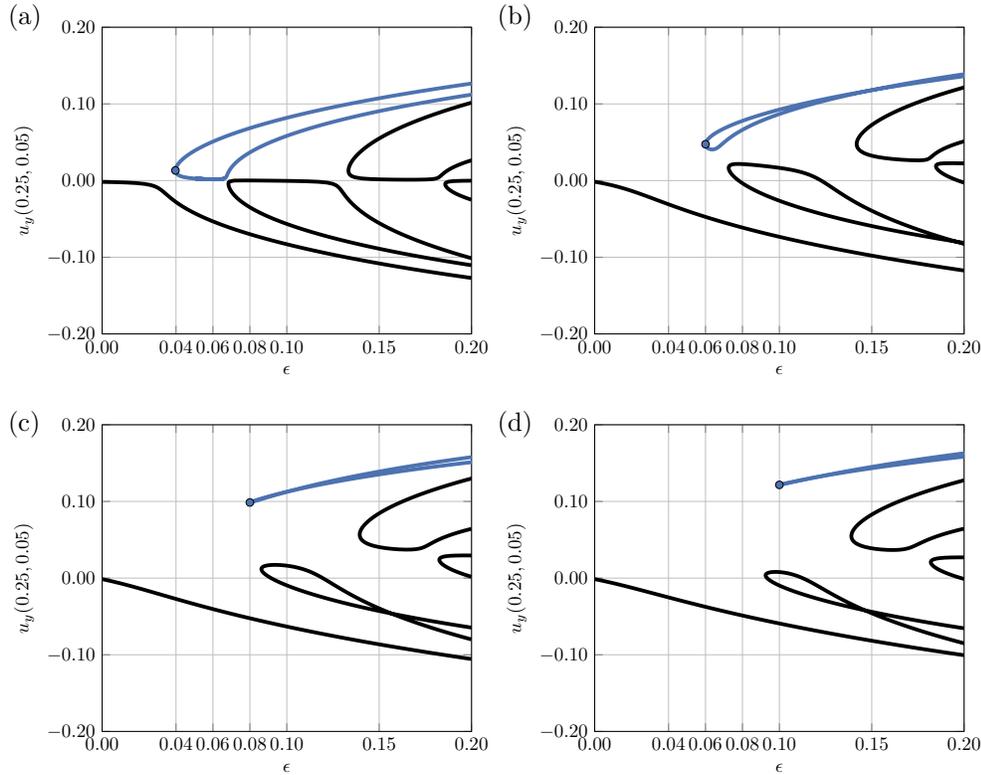
\begin{figure}[htbp]
\vspace{0.2cm}
\begin{center}
\begin{overpic}[width=\textwidth]{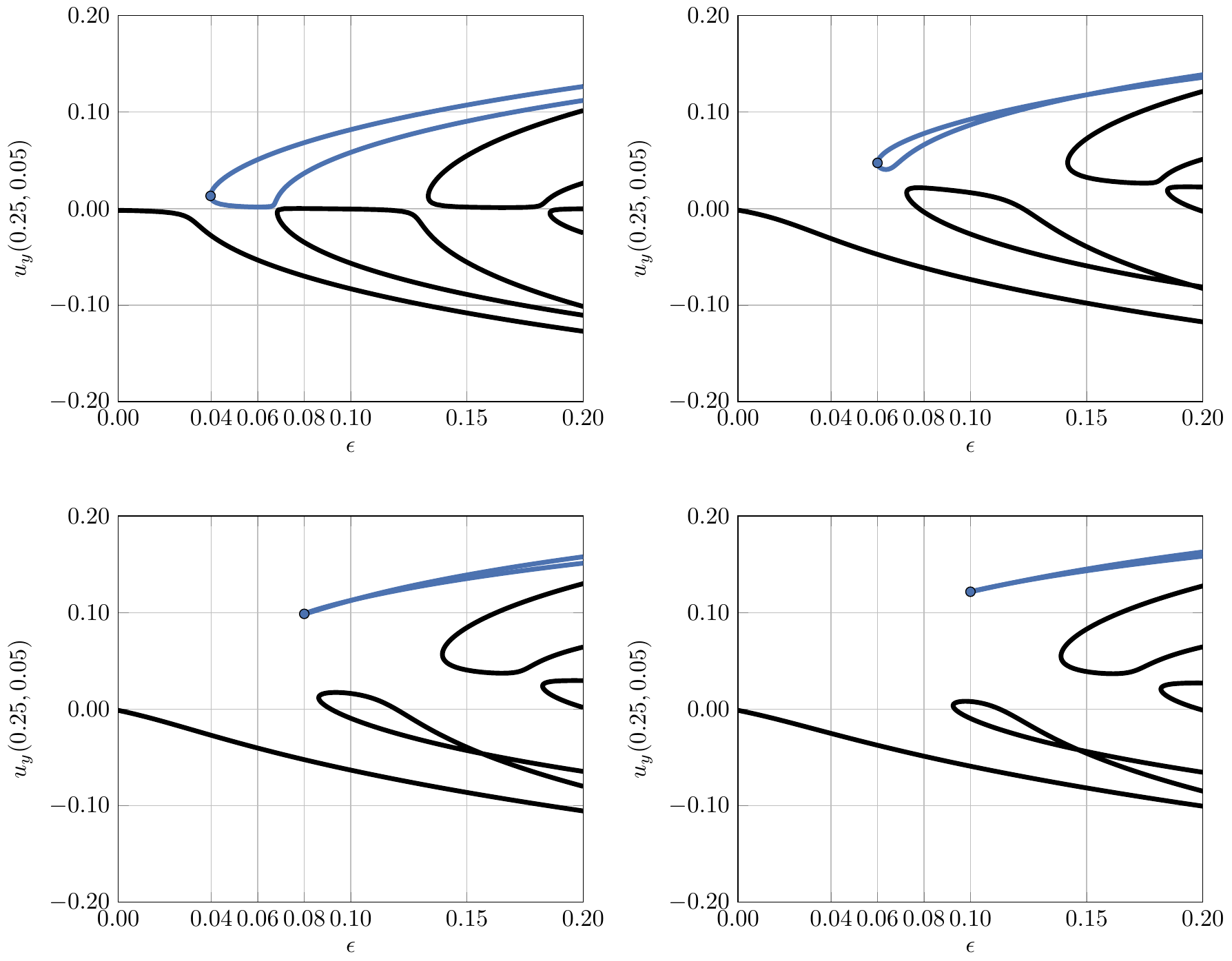}
\put(0,77){(a)}
\put(50,77){(b)}
\put(0,35.5){(c)}
\put(50,35.5){(d)}
\end{overpic}
\end{center}
\caption{Bifurcation diagrams of the hyperelasticity problem computed by deflated continuation, where the first fold bifurcation (in the branch highlighted in blue) arises at $\epsilon\approx 0.04$ (original diagram displayed in panel (a)), $\epsilon=0.06$ (b), $\epsilon=0.08$ (c), and $\epsilon=0.1$ (d). The hyperelasticity equations were solved with the different shapes of beam illustrated in \cref{fig_hyperelasticity_sol}(a)-(d), respectively.}
\label{fig_hyperelasticity}
\end{figure}

The presence of the bifurcation parameter $\epsilon$ in the Dirichlet
boundary conditions presents an additional difficulty. Dirichlet
boundary conditions such as this are typically enforced in the definition of the trial
space, by seeking a solution in
\begin{equation*}
\{v \in H^1(\Omega; \mathbb{R}^2) : \left.v\right|_{\Gamma_D} = g(\epsilon)\},
\end{equation*}
where $g(\epsilon)$ is the given boundary data and $\Gamma_D=\Gamma_L\cup\Gamma_R$
is the boundary on which to apply the conditions.
However, this approach is not feasible when deriving the Moore--Spence
system~\cref{eq_moore_spence} symbolically using the automatic differentiation
facilities of UFL~\cite{alnaes2014}, as these
require having the bifurcation parameter $\epsilon$ as a variable in the weak
form of the equations, rather than in the spaces employed.

To circumvent this difficulty, we choose to seek a solution in
\begin{equation*}
V \coloneqq \{v \in H^1(\Omega; \mathbb{R}^2) : \left.v\right|_{\Gamma_L} = (0, 0)^\top\},
\end{equation*}
and to impose the Dirichlet boundary condition on the
right boundary $\Gamma_R$ using Nitsche's
method~\cite{freund1995weakly,nitsche1971variationsprinzip,sime2020automatic} (as in
Xia et al.~\cite{xia2020nonlinear}). More precisely, we add the following
terms
\begin{equation}\label{eq:Nitsche}
\gamma\int_{\Gamma_R}(u-g(\epsilon))\cdot v\,\mathrm{d}s - \int_{\Gamma_R} (P(u)n)\cdot v\,\mathrm{d}s,
\end{equation}
to the weak formulation arising from the Fr\'echet differentiation of the hyperelastic
energy \cref{eq_energy_hyper}. In formula \cref{eq:Nitsche}, $\gamma>0$ is
the Nitsche penalization parameter and, in the numerical
experiments, we choose
$\gamma=10^{15}$ as in Xia et al.~(justified in~\cite[Appendix~A]{xia2020nonlinear}).
The final weak form is to find $u \in V$ such that
\begin{equation} \label{eq_hyper}
\Pi'(u; v) + \gamma\int_{\Gamma_R}(u-g(\epsilon))\cdot v\,\mathrm{d}s - \int_{\Gamma_R} (P(u)n)\cdot v\,\mathrm{d}s = 0
\end{equation}
for all $v \in V$. We discretize \cref{eq_hyper} using piecewise linear continuous Lagrangian finite elements.

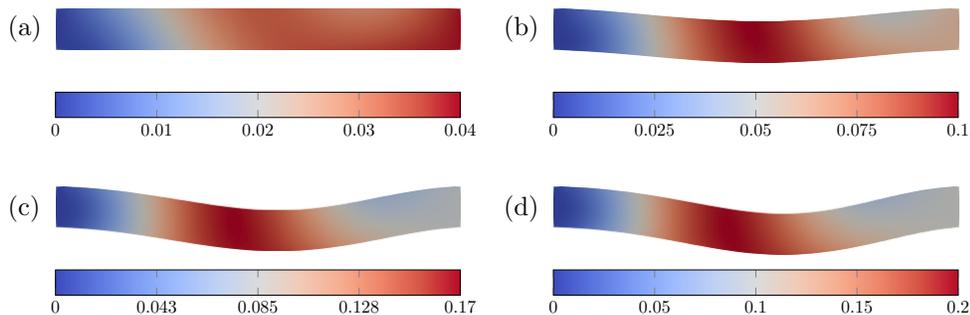
\begin{figure}[htbp]
\vspace{0.5cm}
\begin{center}
\hspace{0.3cm}
\begin{overpic}[width=0.95\textwidth]{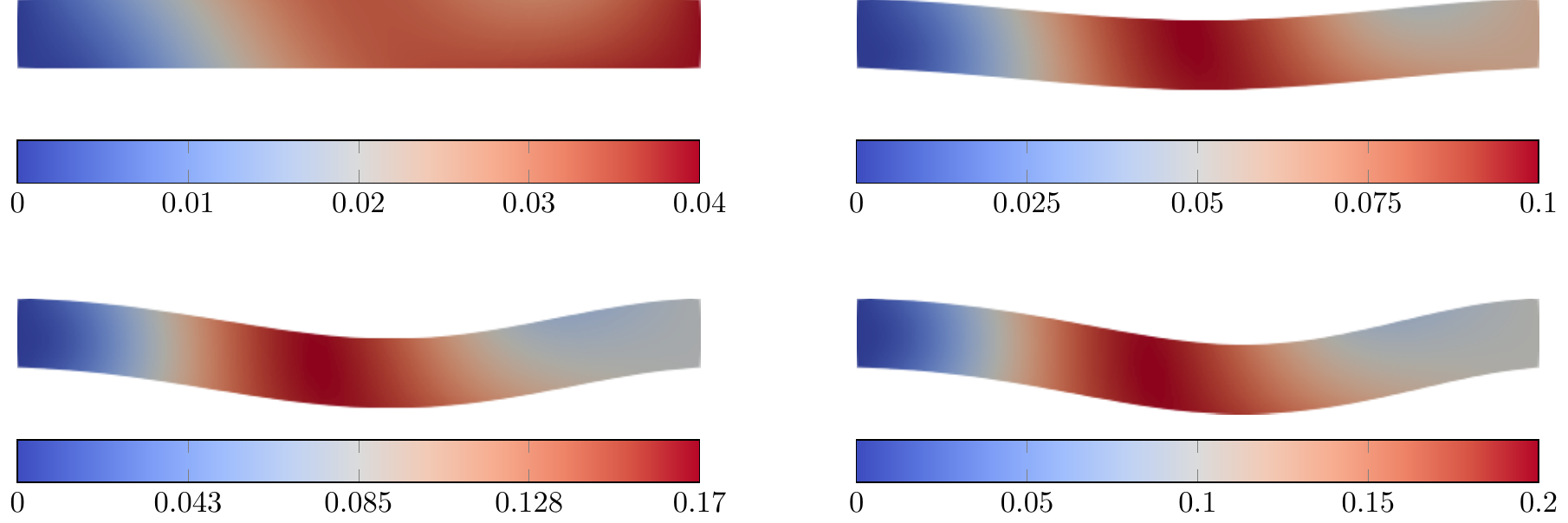}
\put(-4,30.5){(a)}
\put(49.3,30.5){(b)}
\put(-4,11.2){(c)}
\put(49.3,11.2){(d)}
\end{overpic}
\end{center}
\caption{Magnitude of the displacement $u$ at the first fold bifurcation at $\epsilon=0.04$ (a), $\epsilon=0.06$ (b), $\epsilon=0.08$ (c), and $\epsilon=0.1$ (d) on the respective shapes obtained by the optimization algorithm described in \cref{sec_opt_bifurc_points}.}
\label{fig_hyperelasticity_sol}
\end{figure}

We observe several fold
bifurcations as the bifurcation parameter $\epsilon$ varies in the interval
$[0,0.2]$ on the initial domain. The bifurcation diagram of this problem is illustrated in
\cref{fig_hyperelasticity}(a) using the
diagnostic $u_y(0.25,0.05)$ (i.e.~the second component of the displacement) to distinguish the different bifurcation branches. 
In this problem the $\mathbb{Z}_2$-symmetry is broken by the nonzero gravitational body force $B$, causing the structure
to prefer to buckle downwards rather than upwards.
We seek to delay the birth of the
first fold bifurcation arising at $\epsilon\approx 0.04$ to
$\epsilon=0.06,0.08,0.1$ (see~\cref{fig_hyperelasticity}(a) and the
corresponding displacement $u$ illustrated in
\cref{fig_hyperelasticity_sol}(a)). In this example, the left and right boundaries of the domain are fixed during the shape optimization. As in \cref{sec_ns}, the shape optimization algorithm required between 20 and 30 steps 
to minimize the loss function to nearly machine precision in each case.

In \cref{fig_hyperelasticity_sol}(b)-(d), we report the shapes returned by
the shape optimization procedure and the magnitude of the displacement at the
branch points (respectively arising at
$\epsilon=0.06,0.08,0.1$). We observe that the shape of the beam is bent
downward, which makes it difficult to buckle and therefore delays the location of the first fold bifurcation. The bifurcation
diagrams of the hyperelasticity equations solved on these new shapes via
deflated continuation are displayed in \cref{fig_hyperelasticity}(c)-(d). We
find that the bifurcation structure of the problem changes as the location of
the fold bifurcation is delayed and as the shape of the beam is bent. In
\cref{fig_hyperelasticity}(d), we see that the branch that we chose to
control does indeed arise at $\epsilon=0.1$ as desired, but that another fold bifurcation
now exists for a smaller bifurcation parameter. This is because the
numerical technique described in this paper controls a \emph{specific} branch, but
the other branch points are also modified as the shape of the domain changes
during the optimization procedure.

\section{Conclusions}

We have presented a robust computational technique to control the location of a specified branch point using shape optimization. The key idea is to describe the branch point using the Moore--Spence system. We applied this procedure to vary branch points in the bifurcation diagrams of the cubic-quintic Allen--Cahn equation, the incompressible Navier--Stokes equations, and a hyperelasticity equation. Potential future extensions of the algorithm include its generalization to control Hopf bifurcations~\cite{strogatz2018nonlinear} as well as incorporating areas or volume constraints on the domain using augmented Lagrangian optimization algorithms~\cite[Chapt.~17.3]{nocedal2006numerical}. Finally, we expect that this numerical technique could have a wide range of applications in physics and engineering to control the stability of solutions in areas such as the design of new mechanical metamaterials, structural components, flow devices, liquid crystals, and pipes that promote laminar flow.

\section*{Code availability}
For reproducibility, we archived the Firedrake components~\cite{firedrake_zenodo_2021_5217566} and the code~\cite{nicolas_boulle_2021_5235244} used to produce the numerical examples presented in this paper on Zenodo. 
Detailed instructions to reproduce the results are also available on GitHub at \url{https://github.com/NBoulle/Shape_Bifurcation}.

\section*{Acknowledgments}
We thank the referees for their detailed reviews and suggestions.

\bibliographystyle{siam}
\bibliography{biblio}
\end{document}